\documentclass[12pt,a4paper,french,english]{amsart}

\usepackage[utf8]{inputenc}
\usepackage[T1]{fontenc}
\usepackage{babel}
\usepackage{amssymb}
\usepackage{mathrsfs} 
\usepackage{amsmath}
\usepackage{amsthm}
\usepackage{amscd}
\usepackage{pgf,tikz}
\usetikzlibrary{arrows}
\usepackage[alwaysadjust]{paralist} 
\setdefaultenum{\upshape(i)}{}{}{}
\setdefaultleftmargin{10em}{1em}{1em}{1em}{1em}{1em} 
\usepackage[text={16cm,24cm},centering]{geometry}
\usepackage[pdftex,colorlinks,linkcolor=blue,filecolor=red,citecolor=blue,urlcolor=blue]{hyperref} 

\allowdisplaybreaks
\def\oversetalign#1\to#2{\mathbin{\smash{\overset{{\text{\rlap{\hss#1}}}}{#2}}}}
\def\oversettext#1\to#2{\mathbin{\smash{\overset{{\text{{#1}}}}{#2}}}}

\def\dcap_#1{\mathchoice{%
          {\textstyle\bigcap\limits_{#1}}}%
          {\underset{#1}\cap}%
          {\underset{#1}\cap}%
          {\underset{#1}\cap}}
\def\dcup_#1{\mathchoice{%
          {\textstyle\bigcup\limits_{#1}}}%
          {\underset{#1}\cup}%
          {\underset{#1}\cup}%
          {\underset{#1}\cup}}
\def\ddcap_#1^#2{\mathchoice{%
          {\textstyle\bigcap\limits_{#1}^{#2}}}%
          {\overset{#2}{\underset{#1}\cap}}%
          {\overset{#2}{\underset{#1}\cap}}%
          {\overset{#2}{\underset{#1}\cap}}}
\def\ddcup_#1^#2{\mathchoice{%
          {\textstyle\bigcup\limits_{#1}^{#2}}}%
          {\overset{#2}{\underset{#1}\cup}}%
          {\overset{#2}{\underset{#1}\cup}}%
          {\overset{#2}{\underset{#1}\cup}}}

\newskip\myICS\myICS=0.0pt

\newcommand\intervoo[3][1]{\ifcase#1
   \mskip2.0mu\left(\mskip1mu 
   #2,#3
   \mskip1mu\right)\mskip2.0mu
\or 
   \mathopen{\mskip2mu(\mskip1mu}
   #2,#3
   \mathclose{\mskip1mu)\mskip2.0mu}
\or 
   \mathopen{\mskip2mu\bigl(\mskip1mu}
   #2\kern#1\myICS\kern-\myICS,\kern#1\myICS\kern-\myICS#3
   \mathclose{\mskip1mu\bigr)\mskip2.0mu}
\or 
   \mathopen{\mskip2mu\Bigl(\mskip1mu} 
   #2\kern#1\myICS\kern-\myICS,\kern#1\myICS\kern-\myICS#3
   \mathclose{\mskip1mu\Bigr)\mskip2.0mu}
\or 
   \mathopen{\mskip2mu\biggl(\mskip1mu} 
   #2\kern#1\myICS\kern-\myICS,\kern#1\myICS\kern-\myICS#3
   \mathclose{\mskip1mu\biggr)\mskip2.0mu}
\or 
   \mathopen{\mskip2mu\Biggl(\mskip1mu} 
   #2\kern#1\myICS\kern-\myICS,\kern#1\myICS\kern-\myICS#3
   \mathclose{\mskip1mu\Biggr)\mskip2.0mu}
\else 
   \mskip2.0mu\left(\mskip1mu 
   #2\kern#1\myICS\kern-\myICS,\kern#1\myICS\kern-\myICS#3
   \mskip1mu\right)\mskip2.0mu
\fi}

\newcommand\intervof[3][1]{\ifcase#1
   \mskip2.0mu\left(\mskip1mu 
   #2,#3
   \mskip2.5mu\right]\mskip0.5mu
\or 
   \mathopen{\mskip2.0mu(\mskip1mu}
   #2,#3
   \mathclose{\mskip2.5mu]\mskip0.5mu}
\or 
   \mathopen{\mskip2.0mu\bigl(\mskip1mu}
   #2\kern#1\myICS\kern-\myICS,\kern#1\myICS\kern-\myICS#3
   \mathclose{\mskip2.5mu\bigr]\mskip0.5mu}
\or 
   \mathopen{\mskip2mu\Bigl(\mskip1mu} 
   #2\kern#1\myICS\kern-\myICS,\kern#1\myICS\kern-\myICS#3
   \mathclose{\mskip2.5mu\Bigr]\mskip0.5mu}
\or 
   \mathopen{\mskip2.0mu\biggl(\mskip1mu} 
   #2\kern#1\myICS\kern-\myICS,\kern#1\myICS\kern-\myICS#3
   \mathclose{\mskip2.5mu\biggr]\mskip0.5mu}
\or 
   \mathopen{\mskip2.0mu\Biggl(\mskip1mu} 
   #2\kern#1\myICS\kern-\myICS,\kern#1\myICS\kern-\myICS#3
   \mathclose{\mskip2.5mu\Biggr]\mskip0.5mu}
\else 
   \mskip2.0mu\left(\mskip1mu 
   #2\kern#1\myICS\kern-\myICS,\kern#1\myICS\kern-\myICS#3
   \mskip2.5mu\right]\mskip0.5mu
\fi}

\newcommand\intervfo[3][1]{\ifcase#1
   \mskip0.5mu\left[\mskip2mu 
   #2,#3
   \mskip1mu\right)\mskip2.0mu
\or 
   \mathopen{\mskip0.5mu[\mskip2mu}
   #2,#3
   \mathclose{\mskip1mu)\mskip2.0mu}
\or 
   \mathopen{\mskip0.5mu\bigl[\mskip2mu}
   #2\kern#1\myICS\kern-\myICS,\kern#1\myICS\kern-\myICS#3
   \mathclose{\mskip1mu\bigr)\mskip2.0mu}
\or 
   \mathopen{\mskip0.5mu\Bigl[\mskip2mu}
   #2\kern#1\myICS\kern-\myICS,\kern#1\myICS\kern-\myICS#3
   \mathclose{\mskip1mu\Bigr)\mskip2.0mu}
\or 
   \mathopen{\mskip0.5mu\biggl[\mskip2mu}
   #2\kern#1\myICS\kern-\myICS,\kern#1\myICS\kern-\myICS#3
   \mathclose{\mskip1mu\biggr)\mskip2.0mu}
\or 
   \mathopen{\mskip0.5mu\Biggl[\mskip2mu}
   #2\kern#1\myICS\kern-\myICS,\kern#1\myICS\kern-\myICS#3
   \mathclose{\mskip1mu\Biggr)\mskip2.0mu}
\else 
   \mskip0.5mu\left[\mskip2mu 
   #2\kern#1\myICS\kern-\myICS,\kern#1\myICS\kern-\myICS#3
   \mskip1mu\right)\mskip2.0mu
\fi}

\newcommand\intervff[3][1]{\ifcase#1
   \mskip0.5mu\left[\mskip2mu 
   #2,#3
   \mskip2.5mu\right]\mskip0.5mu
\or 
   \mathopen{\mskip0.5mu[\mskip2mu}
   #2,#3
   \mathclose{\mskip2.5mu]\mskip0.5mu}
\or 
   \mathopen{\mskip0.5mu\bigl[\mskip2mu}
   #2\kern#1\myICS\kern-\myICS,\kern#1\myICS\kern-\myICS#3
   \mathclose{\mskip2.5mu\bigr]\mskip0.5mu}
\or 
   \mathopen{\mskip0.5mu\Bigl[\mskip2mu}
   #2\kern#1\myICS\kern-\myICS,\kern#1\myICS\kern-\myICS#3
   \mathclose{\mskip2.5mu\Bigr]\mskip0.5mu}
\or 
   \mathopen{\mskip0.5mu\biggl[\mskip2mu}
   #2\kern#1\myICS\kern-\myICS,\kern#1\myICS\kern-\myICS#3
   \mathclose{\mskip2.5mu\biggr]\mskip0.5mu}
\or 
   \mathopen{\mskip0.5mu\Biggl[\mskip2mu}
   #2\kern#1\myICS\kern-\myICS,\kern#1\myICS\kern-\myICS#3
   \mathclose{\mskip2.5mu\Biggr]\mskip0.5mu}
\else 
   \mskip0.5mu\left[\mskip2mu 
   #2\kern#1\myICS\kern-\myICS,\kern#1\myICS\kern-\myICS#3
   \mskip2.5mu\right]\mskip0.5mu
\fi}

\newdimen\TZschaal\TZschaal=1cm
\def\TZemme{0.2}
\def\TZhoek{25}
\tikzset{param3d/.style={
x={({-sin(\TZhoek)*\TZschaal} , {-\TZemme\TZschaal*cos(\TZhoek)/sqrt(1+\TZemme*\TZemme)})}, 
y={({ cos(\TZhoek)*\TZschaal} , {-\TZemme\TZschaal*sin(\TZhoek)/sqrt(1+\TZemme*\TZemme)})}, 
z={(         0\TZschaal     , {                 \TZschaal/sqrt(1+\TZemme*\TZemme)})}
}}

\newif\ifprintTikZpicture\printTikZpicturetrue
\long\def\WelNietTikZpicture#1{\ifprintTikZpicture#1\fi}

\def\TZVAROBSa{0}\def\TZVAROBSb{3}\def\TZVAROBSc{1}
\def\TZPERSPp#1#2{\ifx h#2 (\TZVAROBSa)+(\TZVAROBSb)*(#1)\else(\TZVAROBSa)+(\TZVAROBSb)*((#1)-(\TZVAROBSa))/(#2)\fi}
\def\TZPERSPq#1#2{\ifx h#2 (\TZVAROBSc)\else(\TZVAROBSc)+(\TZVAROBSb)*((#1)-(\TZVAROBSc))/(#2)\fi}

\newtheoremstyle{mythmstyle}%
{2\baselineskip plus0.15\baselineskip minus 0.15\baselineskip}
{\baselineskip}
{\itshape}
{}
{\bf}
{}
{0pt}
{} 
\newif\ifmynonumberenvi\mynonumberenvitrue
\theoremstyle{mythmstyle}
\newtheorem{proclaimmythm}[equation]{} 
\newtheorem*{proclaimmythm*}{}%

\newenvironment{proclaim}[2][*]{\ifx*#1\mynonumberenvitrue\begin{proclaimmythm*}{\bf#2.} \ignorespaces\else\mynonumberenvifalse\begin{proclaimmythm}{.\kern0.5em\bf#2.}\label{#1} \ignorespaces\fi}{\ifmynonumberenvi\end{proclaimmythm*}\else\end{proclaimmythm}\fi}%

\renewenvironment{proclaim}[2][*]{\ifx*#1\begin{proclaimmythm}{\bf.\kern0.5em\bf#2.} \ignorespaces\else\begin{proclaimmythm}{\bf.\kern0.5em\bf#2.}\label{#1} \ignorespaces\fi}{\end{proclaimmythm}}

\newtheoremstyle{mydefstyle}%
{2\baselineskip plus0.15\baselineskip minus 0.15\baselineskip}
{\baselineskip}
{\rmfamily}
{}
{\bf}
{}
{0pt}
{} 

\theoremstyle{mydefstyle}
\newtheorem{proclaimmydef}[equation]{} 
\newtheorem*{proclaimmydef*}{}%

\newenvironment{definition}[2][*]{\ifx*#1\mynonumberenvitrue\begin{proclaimmydef*}{\bf#2.}\else\mynonumberenvifalse\begin{proclaimmydef}{.\kern0.5em\bf#2.}\label{#1} \ignorespaces\fi}{\ifmynonumberenvi\end{proclaimmydef*}\else\end{proclaimmydef}\fi}

\renewenvironment{definition}[2][*]{\ifx*#1\begin{proclaimmydef}{\bf.\kern0.5em\bf#2.} \ignorespaces\else\begin{proclaimmydef}{\bf.\kern0.5em\bf#2.}\label{#1} \ignorespaces\fi}{\end{proclaimmydef}}

\newenvironment{definition*}[1]{\begin{proclaimmydef*}{\bf#1.} \ignorespaces}{\end{proclaimmydef*}}

\newif\ifSymbolFinRem\SymbolFinRemtrue
\def\FinRemSymbol{\hbox{$\blacktriangleleft\!\blacktriangleright$}}
\def\FinRem{\hskip0.001em plus 40pt\null{}\kern5pt\null\nobreak\hfill
   \kern3pt{\FinRemSymbol}}

\newcommand\FinRemici{\ifSymbolFinRem\\\noalign{\vskip-\baselineskip\smash{\hbox to\linewidth{\vrule width0pt \hfill\global\FinRemdejaplacetrue\FinRemSymbol}}\vskip-\baselineskip}\fi}
\newif\ifFinRemdejaplace\FinRemdejaplacefalse

\newenvironment{remarque}[2][*]%
{\begin{proclaimmydef*}{\bf#2.} \ignorespaces}%
{\ifSymbolFinRem\ifFinRemdejaplace\global\FinRemdejaplacefalse\else\FinRem\fi\fi
\end{proclaimmydef*}}

\newtheoremstyle{mypreuvestyle}%
{\baselineskip}
{\baselineskip}
{}
{}
{\em}
{}
{0pt}
{} 

\def\QEDbox{\hbox{\lower2.3pt\vbox{\hrule\hbox
   {\vrule\kern1pt\vbox{\kern1.7pt\hbox{$\scriptstyle
   QED$}\kern.6pt}\kern1pt\vrule}\hrule}}}
\def\QED{\hskip0.01em plus 40pt\null{} \null\nobreak\hfill
   \kern3pt\QEDbox} 
\newcommand\QEDici{\\\noalign{\vskip-\baselineskip\smash{\hbox to\linewidth{\vrule width0pt \hfill\global\QEDdejaplacetrue\QEDbox}}\vskip-\baselineskip}}
\newif\ifQEDdejaplace\QEDdejaplacefalse

\theoremstyle{mypreuvestyle}
\newtheorem*{proclaimmypreuve}{}
\newenvironment{preuve}[1][*]{\begin{proclaimmypreuve}{\ifx*#1{\em Proof.}\else{\em#1.}\fi} \ignorespaces}{\ifQEDdejaplace\global\QEDdejaplacefalse\else\QED\fi\end{proclaimmypreuve}}

\newif\ifmynonumberequation\mynonumberequationtrue

\makeatletter
\numberwithin{equation}{section}
\newenvironment{moneq}[1][*]{\ifx*#1\mynonumberequationtrue\begin{equation*}\else\mynonumberequationfalse\begin{equation}\label{#1}\fi}{\ifmynonumberequation\end{equation*}\@ignoretrue\else\end{equation}\@ignoretrue\fi\ignorespaces}
\makeatother

\newcommand{\recalf}[1]{(\ref{#1})}
\newcommand{\recals}[1]{{\S\ref{#1}}}
\newcommand{\recalt}[1]{{[\ref{#1}]}}

\newcommand\norme[2][1]{\ifcase#1
   \left\Vert#2\right\Vert 
   \or\Vert#2\Vert 
   \or\bigl\Vert#2\bigr\Vert 
   \or\Bigl\Vert#2\Bigr\Vert 
   \or\biggl\Vert#2\biggr\Vert 
   \or\Biggl\Vert#2\Biggr\Vert 
   \else\Vert#2\Vert\fi} 

\newcommand\inprod[3][1]{\ifcase#1
   \left\langle\mskip1.5mu #2,#3\mskip1.5mu \right\rangle 
   \or\langle\mskip1.5mu #2,#3\mskip1.5mu \rangle 
   \or\bigl\langle\mskip1.5mu #2,#3\mskip1.5mu \bigr\rangle 
   \or\Bigl\langle\mskip1.5mu #2,#3\mskip1.5mu \Bigr\rangle 
   \or\biggl\langle\mskip1.5mu #2,#3\mskip1.5mu \biggr\rangle 
   \or\Biggl\langle\mskip1.5mu #2,#3\mskip1.5mu \Biggr\rangle 
   \else\langle\mskip1.5mu #2,#3\mskip1.5mu \rangle\fi} 

\newcommand\bigrestricted{{\kern1pt\vrule height3.3ex depth1.7ex width0.6pt\kern1pt}}

\newcommand\esptanbas{\operatorname{\mathrm{R}\mkern-4mu\raise0.7ex\hbox{\footnotesize $\mathrm{T}$}\mkern-2mu\mathrm{S}}}
\newcommand\esptanbast{\rlap{\kern0.1em$\widetilde{\kern1.6em\vrule width0pt height1.7ex}$}\esptanbas}
\def\extder{{\operatorname{d}}}
\def\fracp#1#2{\frac{\partial#1}{\partial#2}}

\newcommand\gammah{{\widehat\gamma}}

\newcommand\GL{\mathrm{GL}}
\newcommand\ground{G}
\newcommand\hbiggl[1]{\raise-0.7pt\hbox{\large$\Bigl#1$}}
\newcommand\hbiggr[1]{\raise-0.7pt\hbox{\large$\Bigr#1$}}
\newcommand\hBigl[1]{\raise0.6pt\hbox{\footnotesize$\Bigl#1$}}
\newcommand\hBigr[1]{\raise0.6pt\hbox{\footnotesize$\Bigr#1$}}

\newcommand\ie{{i.e.}}
\newcommand\itemizedef[1][œ]{\ifx#1*\else\par\medskip\noindent\fi$\bullet$ \ignorespaces}
\newcommand\labelPT{PT\kern0.1em\relax}
\newcommand\labelRTS{RTS\kern0.1em\relax}
\newcommand\mapob{\kern0.4em}
\newcommand\Mat{\mathrm{M}}
\newcommand\mfdmetric{\mathbf{g}}
\newcommand\mfdmetrichat{\widehat\mfdmetric}
\def\mo{^{-1}}
\newcommand\myquote[1]{``#1''}
\newcommand\NN{\mathbf{N}}
\newcommand\normal{\mathbf{n}}
\newcommand\refmetnaam[2]{#1\ref{#2}}
\newcommand\RR{\mathbf{R}}
\newcommand\scirc{\,{\raise 0.8pt\hbox{$\scriptstyle\circ$}}\,}
\newcommand\shifttag[1]{\kern#1&\kern-#1}
\newcommand\stresd[1]{{\emph{#1}}}
\newcommand\stress[1]{\textbf{\emph{#1}}}
\newcommand\tcour{t}
\newcommand\tfixe{s}

\newcommand\TraceAff{\mathcal{T}_\tfixe}\newcommand\NoTraceAff[1]{}
\newcommand\transportparallel{\operatorname{\mathrm{P}\mkern-3mu\raise0.7ex\hbox{\footnotesize $\mathrm{T}$}}}
\newcommand\vh{{\widehat v}}

\begin{document}
\renewcommand\labelitemi{\textbullet}

\title[The rolling tangent space]{The rolling tangent space\\ a forgotten vision on parallel transport and geodesics}
\author{Constant Pinteaux}
\email[CP]{FirstName[period]LastName[period]etu[arobase]univ-lille[period]fr}
\author{Gijs M. Tuynman}
\email[GMT]{FirstName[period]LastName[arobase]univ-lille[period]fr}
\address{Univ. Lille, CNRS, UMR 8524 - Laboratoire Paul Painlevé, F-59000 Lille, France
}

\begin{abstract}
Given a submanifold $M\subset \RR^\nu$, a curve $\gamma:I\to M$ and tangent vectors $v$ along $\gamma$, we roll the tangent space along $\gamma$. 
In doing so, we get an imprint\slash trace of $\gamma$ on the tangent space, as well as an imprint\slash trace of the tangent vectors. 
We show that for a vector field $v$ along $\gamma$, the imprint\slash trace of its covariant derivative is the ordinary derivative of its imprint\slash trace vector field. 
It then follows easily that $v$ is a set of parallel vectors along $\gamma$ if and only if their imprint\slash trace on the (affine) tangent space is constant and that $\gamma$ is a geodesic if and only if its trace on the tangent space is a straight line. 

\end{abstract}

\maketitle

\tableofcontents

\section{Introduction}

Despite the fact that its name suggests it is a simple concept, for most students the notion of parallel transport remains mysterious, just as its computational counterpart the covariant derivative (which gets drowned easily in the formalism of connections on (vector) bundles). 
And yet, there exists a way to visualise these notions that makes them easy to understand, a way that was already given in \cite{Persico..1921} and which was taken up by the founder of parallel transport T.~Levi-Civita in his book \cite[Ch. V.10, p101--102]{LeviCivita..1927} on the absolute differential calculus.\footnote{One can find more information on the history of parallel transport in \cite{Tazzioli..2018} or \cite{Tazzioli..2025}; in \cite{CardinTazzioli..2026} one finds interpretations in terms of physics of it.}
Unfortunately, this way to interpret parallel transport seems to be forgotten, as we have not found any trace of it in more \myquote{modern} literature, except in lecture notes of a graduate course by W.T.~van Est given at the University of Amsterdam in the late seventies \cite{VanEst..1979}. 
The purpose of the present text is to put this vision again in the limelight. 

The best way to visualise the idea of the rolling tangent space is the game of a person inside a (big) ball walking on the earth, rolling the ball over the ground. 
When we see the ball\slash submanifold as being fixed, it becomes the tangent space\slash ground that rolls around $M$. 
In this way we obtain a map that starts with a curve $\gamma$ on $M$ and produces a trace curve on the (rolling) tangent space. 
When we have a vector field $v$ along this curve, we also obtain a trace vector field along the trace curve. 
And we show that such a vector field $v$ is \myquote{parallel} if and only if its trace vector field is constant, that the covariant derivative of $v$ is mapped to the ordinary derivative of its trace vector field and, icing on the cake, the curve $\gamma$ is a geodesic if and only if its trace curve in the tangent space is a straight line.

\section{The rolling tangent space}
\label{SECTIONSURPLANTANGENTBASCULANT}

The idea behind the rolling tangent space is the game of a person inside a ball. 
Imagine a big transparent ball (in the form of a rugby ball to make it more interesting), with inside it a person that can walk. 
At each moment this ball touches the ground at a single point: the point where our person inside has its feet. 
And when this person makes a step, he has to \myquote{push} the ball on the ground at a new point, the point where he puts his feet after the step. 
In this way the ball rolls over the ground (without slipping) in the direction taken by the person inside.
\WelNietTikZpicture
{%
\begin{figure}[!ht]
\begin{tikzpicture}[scale=1, line cap=round,line join=round,x=\TZschaal, y=\TZschaal]
\def\TZemme{0.2} 
\def\TZhoek{30} 

\begin{scope}[param3d]

\def\TZVARrectsize{3}
\draw[color=green!70, line width=0.4pt] ({-\TZVARrectsize}, {-\TZVARrectsize}, -1)--({\TZVARrectsize}, {-\TZVARrectsize}, -1)--({\TZVARrectsize}, {\TZVARrectsize}, -1)--({-\TZVARrectsize}, {\TZVARrectsize}, -1)--cycle;

\end{scope}

\draw[color=brown, line width=0.5pt, fill=white] (0,0) circle (1);

\begin{scope}[param3d]

\def\TZVARphione{(-60)}\def\TZVARphitwo{(120)}
\draw[color=brown, line width=0.5pt, smooth, samples=20, domain=\TZVARphione:\TZVARphitwo] plot[variable=\t] ({cos(\t)}, {sin(\t)},{0});
\draw[color=brown, dotted, line width=0.5pt, smooth, samples=20, domain=\TZVARphitwo:\TZVARphione+360] plot[variable=\t] ({cos(\t)}, {sin(\t)},{0});

\draw[color=blue, line width=0.5pt, smooth, samples=20, domain=0:2.5] plot[variable=\t] ({1-cos(2*\t r)-\t},{\t},{-1});

\end{scope}

\begin{scope}[yscale=0.75, scale=0.28]

\def\TZVARpied{0.05} 
\def\TZBHteter{0.12} 
\def\TZBHecartjambes{0.1}
\def\TZBHmiddel{0.45}
\def\TZBHepaule{0.8}
\def\TZBHbras{0.2}\def\TZBHavantbras{0.22}
\def\TZBHcoudegx{\TZBHbras*cos(\TZBHbpg)}\def\TZBHcoudegy{\TZBHbras*sin(\TZBHbpg)}
\def\TZBHmaingx{\TZBHavantbras*cos(\TZBHabpg)}\def\TZBHmaingy{\TZBHavantbras*sin(\TZBHabpg)}
\def\TZBHcoudedx{\TZBHbras*cos(\TZBHbpd)}\def\TZBHcoudedy{\TZBHbras*sin(\TZBHbpd)}
\def\TZBHmaindx{\TZBHavantbras*cos(\TZBHabpd)}\def\TZBHmaindy{\TZBHavantbras*sin(\TZBHabpd)}
\def\TZBHabpg{-80}\def\TZBHbpg{20}\def\TZBHabpd{110}\def\TZBHbpd{30}

\def\TZVARoeilz{1}
\def\TZVARoeily{0.55}
\def\TZVARoeilx{0}
\def\TZBHabpg{-80}\def\TZBHbpg{20}\def\TZBHabpd{130}\def\TZBHbpd{30}

\draw ({\TZPERSPp{\TZVARoeilx}{\TZVARoeily}},{\TZPERSPq{\TZBHmiddel}{\TZVARoeily}})--({\TZPERSPp{\TZVARoeilx}{\TZVARoeily}},{\TZPERSPq{\TZVARoeilz}{\TZVARoeily}});

\draw ({\TZPERSPp{\TZVARoeilx-\TZBHecartjambes-\TZVARpied}{\TZVARoeily}},{\TZPERSPq{0}{\TZVARoeily}})--({\TZPERSPp{\TZVARoeilx-\TZBHecartjambes}{\TZVARoeily}},{\TZPERSPq{0}{\TZVARoeily}})--({\TZPERSPp{\TZVARoeilx}{\TZVARoeily}},{\TZPERSPq{\TZBHmiddel}{\TZVARoeily}})--({\TZPERSPp{\TZVARoeilx+\TZBHecartjambes}{\TZVARoeily}},{\TZPERSPq{0}{\TZVARoeily}})--({\TZPERSPp{\TZVARoeilx+\TZBHecartjambes-\TZVARpied}{\TZVARoeily}},{\TZPERSPq{0}{\TZVARoeily}});

\draw 
({\TZPERSPp{\TZVARoeilx-\TZBHmaingx-\TZBHcoudegx}{\TZVARoeily}},{\TZPERSPq{\TZBHepaule-\TZBHmaingy-\TZBHcoudegy}{\TZVARoeily}})
--
({\TZPERSPp{\TZVARoeilx-\TZBHcoudegx}{\TZVARoeily}},{\TZPERSPq{\TZBHepaule-\TZBHcoudegy}{\TZVARoeily}})
--
({\TZPERSPp{\TZVARoeilx}{\TZVARoeily}},{\TZPERSPq{\TZBHepaule}{\TZVARoeily}})
--
({\TZPERSPp{\TZVARoeilx+\TZBHcoudedx}{\TZVARoeily}},{\TZPERSPq{\TZBHepaule-\TZBHcoudedy}{\TZVARoeily}})
--
({\TZPERSPp{\TZVARoeilx+\TZBHmaindx+\TZBHcoudedx}{\TZVARoeily}},{\TZPERSPq{\TZBHepaule-\TZBHmaindy-\TZBHcoudedy}{\TZVARoeily}})
;

\def\TZVARx{\TZVARoeilx+\TZBHteter*cos(\t)}
\def\TZVARy{\TZVARoeily}
\def\TZVARz{\TZVARoeilz+\TZBHteter*sin(\t)}
\def\TZVARborne{360}
\draw[fill=lightgray, line width=0.5pt, smooth, samples=30, domain=0:\TZVARborne] plot[variable=\t] ({\TZPERSPp{\TZVARx}{\TZVARy}}, {\TZPERSPq{\TZVARz}{\TZVARy}});

\end{scope} 

\end{tikzpicture}
\end{figure}
}%
Now suppose we have painted a curve on this (transparent) ball and suppose the person inside follows this curve: at each step he puts his feet a bit farther on this curve. 
In this way he will follow this curve, the ball will roll over te ground and he will trace a new curve on the ground: the points where his feet are at every step. 
For instance, if the ball is a (perfect) sphere and the painted curve is a great circle (equator), then this trace on the ground will be a straight line. 
But if the painted curve is a (small) circle, then the trace curve on the ground will also be a circle. 
Another example is a cone (Chinese hat). 
If the painted curve is a circle around the tip of the cone, then (as for a small circle on a sphere) the trace curve on the ground will also be a circle. 
And the movement of the cone on the ground is described by a rotation around the tip of the cone. 
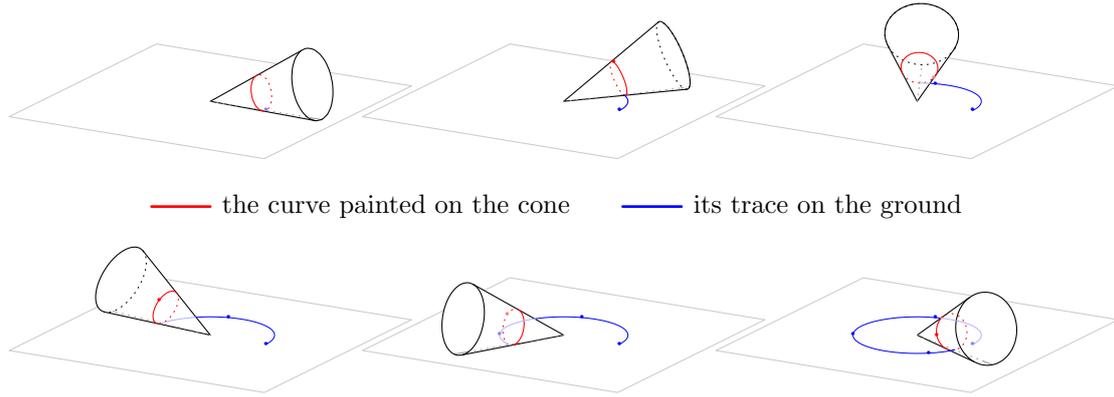
\WelNietTikZpicture
{%
\begin{figure}[!ht]
\begin{tikzpicture}[scale=1.55, line cap=round,line join=round,x=\TZschaal, y=\TZschaal]

\def\TZemme{0.3}
\def\TZhoek{30}

\def\TZVARbasesquare{1.25}

\def\TZVARh{0.3}
\def\TZVARhsqrt{(\TZVARh*sqrt(1+(\TZVARh)*(\TZVARh)))}

\def\TZFx#1#2#3{%
(#1)*(cos(#2 r) - \TZVARhsqrt*sin(#3 r)*sin(#2 r)+(\TZVARh)*(\TZVARh)*cos(#3 r)*cos(#2 r))
}
\def\TZFy#1#2#3{%
(#1)*(sin(#2 r) + \TZVARhsqrt*sin(#3 r)*cos(#2 r)+(\TZVARh)*(\TZVARh)*cos(#3 r)*sin(#2 r))
}
\def\TZFz#1#2#3{%
(#1)*(\TZVARh)*(1-cos(#3 r))
}

\begin{scope}

\def\TZVARrho{1}
\def\TZVARphi{pi/2}
\def\TZVARpsi{0}

\begin{scope}[param3d]

\draw[color=black!20] (-\TZVARbasesquare,-\TZVARbasesquare,0)--(\TZVARbasesquare,-\TZVARbasesquare,0)--(\TZVARbasesquare,\TZVARbasesquare,0)--(-\TZVARbasesquare,\TZVARbasesquare,0)--cycle;

\def\TZVARrho{1}
\def\TZVARrhointer{0.5}
\def\TZVARphi{(pi/2)}
\def\TZVARpsi{0}

\def\TZVARpsibordone{3.25}
\def\TZVARpsibordtwo{-0.3}
\fill[color=white] plot[variable=\TZVARpsi, smooth, samples=20, domain=0:2*pi] ({\TZFx{\TZVARrho}{\TZVARphi}{\TZVARpsi}},{\TZFy{\TZVARrho}{\TZVARphi}{\TZVARpsi}},{\TZFz{\TZVARrho}{\TZVARphi}{\TZVARpsi}});

\draw[color=black, fill=white, line width=0.4pt] ({\TZFx{\TZVARrho}{\TZVARphi}{\TZVARpsibordone}},{\TZFy{\TZVARrho}{\TZVARphi}{\TZVARpsibordone}},{\TZFz{\TZVARrho}{\TZVARphi}{\TZVARpsibordone}})--(0,0,0)--({\TZFx{\TZVARrho}{\TZVARphi}{\TZVARpsibordtwo}},{\TZFy{\TZVARrho}{\TZVARphi}{\TZVARpsibordtwo}},{\TZFz{\TZVARrho}{\TZVARphi}{\TZVARpsibordtwo}});

\def\TZVARpsi{0}
\def\TZVARrhointer{0.9}
\draw[color=black!40, dotted, line width=0.4pt] (0,0,0)--({\TZFx{\TZVARrhointer}{\TZVARphi}{\TZVARpsi}},{\TZFy{\TZVARrhointer}{\TZVARphi}{\TZVARpsi}},{\TZFz{\TZVARrhointer}{\TZVARphi}{\TZVARpsi}});
\draw[color=black!40, line width=0.4pt] ({\TZFx{\TZVARrhointer}{\TZVARphi}{\TZVARpsi}},{\TZFy{\TZVARrhointer}{\TZVARphi}{\TZVARpsi}},{\TZFz{\TZVARrhointer}{\TZVARphi}{\TZVARpsi}})--({\TZFx{\TZVARrho}{\TZVARphi}{\TZVARpsi}},{\TZFy{\TZVARrho}{\TZVARphi}{\TZVARpsi}},{\TZFz{\TZVARrho}{\TZVARphi}{\TZVARpsi}});

\foreach\TZVARrho in {1}
\draw[color=black, line width=0.4pt, smooth, samples=20, domain=0:2*pi] plot[variable=\TZVARpsi] ({\TZFx{\TZVARrho}{\TZVARphi}{\TZVARpsi}},{\TZFy{\TZVARrho}{\TZVARphi}{\TZVARpsi}},{\TZFz{\TZVARrho}{\TZVARphi}{\TZVARpsi}});

\def\TZVARpsibordone{-0.35}
\def\TZVARpsibordtwo{3.1}
\foreach\TZVARrho in { 0.5}
{
\draw[color=red, line width=0.4pt, dotted, smooth, samples=20, domain=\TZVARpsibordone:\TZVARpsibordtwo] plot[variable=\TZVARpsi] ({\TZFx{\TZVARrho}{\TZVARphi}{\TZVARpsi}},{\TZFy{\TZVARrho}{\TZVARphi}{\TZVARpsi}},{\TZFz{\TZVARrho}{\TZVARphi}{\TZVARpsi}});
\draw[color=red, line width=0.4pt, smooth, samples=20, domain=\TZVARpsibordtwo:\TZVARpsibordone+2*pi] plot[variable=\TZVARpsi] ({\TZFx{\TZVARrho}{\TZVARphi}{\TZVARpsi}},{\TZFy{\TZVARrho}{\TZVARphi}{\TZVARpsi}},{\TZFz{\TZVARrho}{\TZVARphi}{\TZVARpsi}});
}

\def\TZVARrhointer{0.5}\def\TZVARpsi{0}

\draw[color=blue!50] ({\TZFx{\TZVARrhointer}{\TZVARphi}{\TZVARpsi}},{\TZFy{\TZVARrhointer}{\TZVARphi}{\TZVARpsi}},{\TZFz{\TZVARrhointer}{\TZVARphi}{\TZVARpsi}}) node {$\cdot$};

\end{scope}

\end{scope}

\begin{scope}[xshift=3cm]

\def\TZVARrho{1}
\def\TZVARphi{pi/2}
\def\TZVARpsi{0}

\begin{scope}[param3d]

\draw[color=black!20] (-\TZVARbasesquare,-\TZVARbasesquare,0)--(\TZVARbasesquare,-\TZVARbasesquare,0)--(\TZVARbasesquare,\TZVARbasesquare,0)--(-\TZVARbasesquare,\TZVARbasesquare,0)--cycle;

\def\TZVARphi{(pi/2+0.9)}
\def\TZVARpsi{0}

\fill[color=white] plot[variable=\TZVARpsi, smooth, samples=20, domain=0:2*pi] ({\TZFx{\TZVARrho}{\TZVARphi}{\TZVARpsi}},{\TZFy{\TZVARrho}{\TZVARphi}{\TZVARpsi}},{\TZFz{\TZVARrho}{\TZVARphi}{\TZVARpsi}});

\def\TZVARpsibordone{2.6}
\def\TZVARpsibordtwo{-0.3}
\draw[color=black, fill=white, line width=0.4pt] ({\TZFx{\TZVARrho}{\TZVARphi}{\TZVARpsibordone}},{\TZFy{\TZVARrho}{\TZVARphi}{\TZVARpsibordone}},{\TZFz{\TZVARrho}{\TZVARphi}{\TZVARpsibordone}})--(0,0,0)--({\TZFx{\TZVARrho}{\TZVARphi}{\TZVARpsibordtwo}},{\TZFy{\TZVARrho}{\TZVARphi}{\TZVARpsibordtwo}},{\TZFz{\TZVARrho}{\TZVARphi}{\TZVARpsibordtwo}});

\def\TZVARpsi{0}
\def\TZVARrhointer{0.98}
\draw[color=black!40, dotted, line width=0.4pt] (0,0,0)--({\TZFx{\TZVARrhointer}{\TZVARphi}{\TZVARpsi}},{\TZFy{\TZVARrhointer}{\TZVARphi}{\TZVARpsi}},{\TZFz{\TZVARrhointer}{\TZVARphi}{\TZVARpsi}});
\draw[color=black!40, line width=0.4pt] ({\TZFx{\TZVARrhointer}{\TZVARphi}{\TZVARpsi}},{\TZFy{\TZVARrhointer}{\TZVARphi}{\TZVARpsi}},{\TZFz{\TZVARrhointer}{\TZVARphi}{\TZVARpsi}})--({\TZFx{\TZVARrho}{\TZVARphi}{\TZVARpsi}},{\TZFy{\TZVARrho}{\TZVARphi}{\TZVARpsi}},{\TZFz{\TZVARrho}{\TZVARphi}{\TZVARpsi}});

\foreach\TZVARrho in { 1}
{
\draw[color=black, line width=0.4pt, dotted, smooth, samples=20, domain=\TZVARpsibordtwo:\TZVARpsibordtwo+2*pi] plot[variable=\TZVARpsi] ({\TZFx{\TZVARrho}{\TZVARphi}{\TZVARpsi}},{\TZFy{\TZVARrho}{\TZVARphi}{\TZVARpsi}},{\TZFz{\TZVARrho}{\TZVARphi}{\TZVARpsi}});
\draw[color=black, line width=0.4pt, smooth, samples=20, domain=\TZVARpsibordone:\TZVARpsibordtwo+2*pi] plot[variable=\TZVARpsi] ({\TZFx{\TZVARrho}{\TZVARphi}{\TZVARpsi}},{\TZFy{\TZVARrho}{\TZVARphi}{\TZVARpsi}},{\TZFz{\TZVARrho}{\TZVARphi}{\TZVARpsi}});
}

\def\TZVARpsibordone{-0.3}
\def\TZVARpsibordtwo{2.6}
\foreach\TZVARrho in { 0.5}
{
\draw[color=red, line width=0.4pt, dotted, smooth, samples=20, domain=\TZVARpsibordone:\TZVARpsibordtwo] plot[variable=\TZVARpsi] ({\TZFx{\TZVARrho}{\TZVARphi}{\TZVARpsi}},{\TZFy{\TZVARrho}{\TZVARphi}{\TZVARpsi}},{\TZFz{\TZVARrho}{\TZVARphi}{\TZVARpsi}});
\draw[color=red, line width=0.4pt, smooth, samples=20, domain=\TZVARpsibordtwo:\TZVARpsibordone+2*pi] plot[variable=\TZVARpsi] ({\TZFx{\TZVARrho}{\TZVARphi}{\TZVARpsi}},{\TZFy{\TZVARrho}{\TZVARphi}{\TZVARpsi}},{\TZFz{\TZVARrho}{\TZVARphi}{\TZVARpsi}});
}

\def\TZVARrhointer{0.5}
\draw[color=blue, line width=0.4pt, smooth, samples=20, domain=(pi/2):\TZVARphi] plot[variable=\TZVARphi] ({\TZFx{\TZVARrhointer}{\TZVARphi}{\TZVARpsi}},{\TZFy{\TZVARrhointer}{\TZVARphi}{\TZVARpsi}},{\TZFz{\TZVARrhointer}{\TZVARphi}{\TZVARpsi}});

\def\TZVARpsi{(((pi/2)-\TZVARphi)*sqrt(1+\TZVARh*\TZVARh)/\TZVARh)}
\draw[color=red] ({\TZFx{\TZVARrhointer}{\TZVARphi}{\TZVARpsi}},{\TZFy{\TZVARrhointer}{\TZVARphi}{\TZVARpsi}},{\TZFz{\TZVARrhointer}{\TZVARphi}{\TZVARpsi}}) node {$\cdot$};

\def\TZVARpsi{0}\def\TZVARphi{pi/2}
\draw[color=blue] ({\TZFx{\TZVARrhointer}{\TZVARphi}{\TZVARpsi}},{\TZFy{\TZVARrhointer}{\TZVARphi}{\TZVARpsi}},{\TZFz{\TZVARrhointer}{\TZVARphi}{\TZVARpsi}}) node {$\cdot$};

\end{scope}

\end{scope}

\begin{scope}[xshift=6cm]

\def\TZVARrho{1}
\def\TZVARphi{pi/2}
\def\TZVARpsi{0}

\begin{scope}[param3d]

\draw[color=black!20] (-\TZVARbasesquare,-\TZVARbasesquare,0)--(\TZVARbasesquare,-\TZVARbasesquare,0)--(\TZVARbasesquare,\TZVARbasesquare,0)--(-\TZVARbasesquare,\TZVARbasesquare,0)--cycle;

\def\TZVARphi{(5*pi/4-0.3)}
\def\TZVARpsi{0}

\fill[color=white] plot[variable=\TZVARpsi, smooth, samples=20, domain=0:2*pi] ({\TZFx{\TZVARrho}{\TZVARphi}{\TZVARpsi}},{\TZFy{\TZVARrho}{\TZVARphi}{\TZVARpsi}},{\TZFz{\TZVARrho}{\TZVARphi}{\TZVARpsi}});

\def\TZVARpsibordone{1.2}
\def\TZVARpsibordtwo{-1.1}
\draw[color=black, fill=white, line width=0.4pt] ({\TZFx{\TZVARrho}{\TZVARphi}{\TZVARpsibordone}},{\TZFy{\TZVARrho}{\TZVARphi}{\TZVARpsibordone}},{\TZFz{\TZVARrho}{\TZVARphi}{\TZVARpsibordone}})--(0,0,0)--({\TZFx{\TZVARrho}{\TZVARphi}{\TZVARpsibordtwo}},{\TZFy{\TZVARrho}{\TZVARphi}{\TZVARpsibordtwo}},{\TZFz{\TZVARrho}{\TZVARphi}{\TZVARpsibordtwo}});

\def\TZVARpsi{0}
\def\TZVARrhointer{0.98}
\draw[color=black!40, dotted, line width=0.4pt] (0,0,0)--({\TZFx{\TZVARrhointer}{\TZVARphi}{\TZVARpsi}},{\TZFy{\TZVARrhointer}{\TZVARphi}{\TZVARpsi}},{\TZFz{\TZVARrhointer}{\TZVARphi}{\TZVARpsi}});
\draw[color=black!40, line width=0.4pt] ({\TZFx{\TZVARrhointer}{\TZVARphi}{\TZVARpsi}},{\TZFy{\TZVARrhointer}{\TZVARphi}{\TZVARpsi}},{\TZFz{\TZVARrhointer}{\TZVARphi}{\TZVARpsi}})--({\TZFx{\TZVARrho}{\TZVARphi}{\TZVARpsi}},{\TZFy{\TZVARrho}{\TZVARphi}{\TZVARpsi}},{\TZFz{\TZVARrho}{\TZVARphi}{\TZVARpsi}});

\foreach\TZVARrho in { 1}
{
\draw[color=black, line width=0.4pt, dotted, smooth, samples=20, domain=\TZVARpsibordtwo:\TZVARpsibordtwo+2*pi] plot[variable=\TZVARpsi] ({\TZFx{\TZVARrho}{\TZVARphi}{\TZVARpsi}},{\TZFy{\TZVARrho}{\TZVARphi}{\TZVARpsi}},{\TZFz{\TZVARrho}{\TZVARphi}{\TZVARpsi}});
\draw[color=black, line width=0.4pt, smooth, samples=20, domain=\TZVARpsibordone:\TZVARpsibordtwo+2*pi] plot[variable=\TZVARpsi] ({\TZFx{\TZVARrho}{\TZVARphi}{\TZVARpsi}},{\TZFy{\TZVARrho}{\TZVARphi}{\TZVARpsi}},{\TZFz{\TZVARrho}{\TZVARphi}{\TZVARpsi}});
}

\def\TZVARpsibordone{-1.1}
\def\TZVARpsibordtwo{1.2}
\foreach\TZVARrho in { 0.5}
{
\draw[color=red, line width=0.4pt, dotted, smooth, samples=20, domain=\TZVARpsibordone:\TZVARpsibordtwo] plot[variable=\TZVARpsi] ({\TZFx{\TZVARrho}{\TZVARphi}{\TZVARpsi}},{\TZFy{\TZVARrho}{\TZVARphi}{\TZVARpsi}},{\TZFz{\TZVARrho}{\TZVARphi}{\TZVARpsi}});
\draw[color=red, line width=0.4pt, smooth, samples=20, domain=\TZVARpsibordtwo:\TZVARpsibordone+2*pi] plot[variable=\TZVARpsi] ({\TZFx{\TZVARrho}{\TZVARphi}{\TZVARpsi}},{\TZFy{\TZVARrho}{\TZVARphi}{\TZVARpsi}},{\TZFz{\TZVARrho}{\TZVARphi}{\TZVARpsi}});
}

\def\TZVARrhointer{0.5}
\def\TZVARphiinter{3.46}
\draw[color=blue, line width=0.3pt, smooth, samples=20, domain=(pi/2):\TZVARphiinter] plot[variable=\TZVARphi] ({\TZFx{\TZVARrhointer}{\TZVARphi}{\TZVARpsi}},{\TZFy{\TZVARrhointer}{\TZVARphi}{\TZVARpsi}},{\TZFz{\TZVARrhointer}{\TZVARphi}{\TZVARpsi}});
\draw[color=blue!30, line width=0.3pt, smooth, samples=20, domain=\TZVARphiinter:\TZVARphi] plot[variable=\TZVARphi] ({\TZFx{\TZVARrhointer}{\TZVARphi}{\TZVARpsi}},{\TZFy{\TZVARrhointer}{\TZVARphi}{\TZVARpsi}},{\TZFz{\TZVARrhointer}{\TZVARphi}{\TZVARpsi}});

\def\TZVARpsi{(((pi/2)-\TZVARphi)*sqrt(1+\TZVARh*\TZVARh)/\TZVARh)}
\draw[color=red!50] ({\TZFx{\TZVARrhointer}{\TZVARphi}{\TZVARpsi}},{\TZFy{\TZVARrhointer}{\TZVARphi}{\TZVARpsi}},{\TZFz{\TZVARrhointer}{\TZVARphi}{\TZVARpsi}}) node {$\cdot$};

\def\TZVARpsi{0}\def\TZVARphi{pi/2}
\draw[color=blue] ({\TZFx{\TZVARrhointer}{\TZVARphi}{\TZVARpsi}},{\TZFy{\TZVARrhointer}{\TZVARphi}{\TZVARpsi}},{\TZFz{\TZVARrhointer}{\TZVARphi}{\TZVARpsi}}) node {$\cdot$};

\def\TZVARpsi{0}\def\TZVARphi{((pi/2)+2*pi*\TZVARh/sqrt(1+\TZVARh*\TZVARh))}
\draw[color=blue] ({\TZFx{\TZVARrhointer}{\TZVARphi}{\TZVARpsi}},{\TZFy{\TZVARrhointer}{\TZVARphi}{\TZVARpsi}},{\TZFz{\TZVARrhointer}{\TZVARphi}{\TZVARpsi}}) node {$\cdot$};

\end{scope}

\end{scope}

\def\TZVARdescente{-0.9}
\draw[color=red, line width=1pt] (-0.5,\TZVARdescente)--++(0.5,0) node[right, color=black] {\footnotesize the curve on the cone};
\begin{scope}[xshift=4cm]
\draw[color=blue, line width=1pt] (-0.5,\TZVARdescente)--++(0.5,0) node[right, color=black] {\footnotesize the trace on the ground};
\end{scope}

\begin{scope}[yshift=-2cm, xshift=0cm]

\def\TZVARrho{1}
\def\TZVARphi{pi/2}
\def\TZVARpsi{0}

\begin{scope}[param3d]

\draw[color=black!20] (-\TZVARbasesquare,-\TZVARbasesquare,0)--(\TZVARbasesquare,-\TZVARbasesquare,0)--(\TZVARbasesquare,\TZVARbasesquare,0)--(-\TZVARbasesquare,\TZVARbasesquare,0)--cycle;

\def\TZVARphi{(5*pi/4+0.6)}
\def\TZVARpsi{0}

\fill[color=white] plot[variable=\TZVARpsi, smooth, samples=20, domain=0:2*pi] ({\TZFx{\TZVARrho}{\TZVARphi}{\TZVARpsi}},{\TZFy{\TZVARrho}{\TZVARphi}{\TZVARpsi}},{\TZFz{\TZVARrho}{\TZVARphi}{\TZVARpsi}});

\def\TZVARpsibordone{-2.4}
\def\TZVARpsibordtwo{0.4}
\draw[color=black, fill=white, line width=0.4pt] ({\TZFx{\TZVARrho}{\TZVARphi}{\TZVARpsibordone}},{\TZFy{\TZVARrho}{\TZVARphi}{\TZVARpsibordone}},{\TZFz{\TZVARrho}{\TZVARphi}{\TZVARpsibordone}})--(0,0,0)--({\TZFx{\TZVARrho}{\TZVARphi}{\TZVARpsibordtwo}},{\TZFy{\TZVARrho}{\TZVARphi}{\TZVARpsibordtwo}},{\TZFz{\TZVARrho}{\TZVARphi}{\TZVARpsibordtwo}});

\def\TZVARpsi{0}
\def\TZVARrhointer{0.98}
\draw[color=black!40, dotted, line width=0.4pt] (0,0,0)--({\TZFx{\TZVARrhointer}{\TZVARphi}{\TZVARpsi}},{\TZFy{\TZVARrhointer}{\TZVARphi}{\TZVARpsi}},{\TZFz{\TZVARrhointer}{\TZVARphi}{\TZVARpsi}});
\draw[color=black!40, line width=0.4pt] ({\TZFx{\TZVARrhointer}{\TZVARphi}{\TZVARpsi}},{\TZFy{\TZVARrhointer}{\TZVARphi}{\TZVARpsi}},{\TZFz{\TZVARrhointer}{\TZVARphi}{\TZVARpsi}})--({\TZFx{\TZVARrho}{\TZVARphi}{\TZVARpsi}},{\TZFy{\TZVARrho}{\TZVARphi}{\TZVARpsi}},{\TZFz{\TZVARrho}{\TZVARphi}{\TZVARpsi}});

\foreach\TZVARrho in { 1}
{
\draw[color=black, line width=0.4pt, dotted, smooth, samples=20, domain=\TZVARpsibordone:\TZVARpsibordtwo] plot[variable=\TZVARpsi] ({\TZFx{\TZVARrho}{\TZVARphi}{\TZVARpsi}},{\TZFy{\TZVARrho}{\TZVARphi}{\TZVARpsi}},{\TZFz{\TZVARrho}{\TZVARphi}{\TZVARpsi}});
\draw[color=black, line width=0.4pt, smooth, samples=20, domain=\TZVARpsibordtwo:\TZVARpsibordone+2*pi] plot[variable=\TZVARpsi] ({\TZFx{\TZVARrho}{\TZVARphi}{\TZVARpsi}},{\TZFy{\TZVARrho}{\TZVARphi}{\TZVARpsi}},{\TZFz{\TZVARrho}{\TZVARphi}{\TZVARpsi}});
}

\def\TZVARpsibordone{-2.4}
\def\TZVARpsibordtwo{0.4}
\foreach\TZVARrho in { 0.5}
{
\draw[color=red, line width=0.4pt, dotted, smooth, samples=20, domain=\TZVARpsibordone:\TZVARpsibordtwo] plot[variable=\TZVARpsi] ({\TZFx{\TZVARrho}{\TZVARphi}{\TZVARpsi}},{\TZFy{\TZVARrho}{\TZVARphi}{\TZVARpsi}},{\TZFz{\TZVARrho}{\TZVARphi}{\TZVARpsi}});
\draw[color=red, line width=0.4pt, smooth, samples=20, domain=\TZVARpsibordtwo:\TZVARpsibordone+2*pi] plot[variable=\TZVARpsi] ({\TZFx{\TZVARrho}{\TZVARphi}{\TZVARpsi}},{\TZFy{\TZVARrho}{\TZVARphi}{\TZVARpsi}},{\TZFz{\TZVARrho}{\TZVARphi}{\TZVARpsi}});
}

\def\TZVARrhointer{0.5}
\def\TZVARphiinter{3.9}
\draw[color=blue, line width=0.3pt, smooth, samples=20, domain=(pi/2):\TZVARphiinter] plot[variable=\TZVARphi] ({\TZFx{\TZVARrhointer}{\TZVARphi}{\TZVARpsi}},{\TZFy{\TZVARrhointer}{\TZVARphi}{\TZVARpsi}},{\TZFz{\TZVARrhointer}{\TZVARphi}{\TZVARpsi}});
\draw[color=blue!30, line width=0.2pt, smooth, samples=20, domain=\TZVARphiinter:\TZVARphi] plot[variable=\TZVARphi] ({\TZFx{\TZVARrhointer}{\TZVARphi}{\TZVARpsi}},{\TZFy{\TZVARrhointer}{\TZVARphi}{\TZVARpsi}},{\TZFz{\TZVARrhointer}{\TZVARphi}{\TZVARpsi}});

\def\TZVARpsi{(((pi/2)-\TZVARphi)*sqrt(1+\TZVARh*\TZVARh)/\TZVARh)}
\draw[color=red] ({\TZFx{\TZVARrhointer}{\TZVARphi}{\TZVARpsi}},{\TZFy{\TZVARrhointer}{\TZVARphi}{\TZVARpsi}},{\TZFz{\TZVARrhointer}{\TZVARphi}{\TZVARpsi}}) node {$\cdot$};

\def\TZVARpsi{0}\def\TZVARphi{pi/2}
\draw[color=blue] ({\TZFx{\TZVARrhointer}{\TZVARphi}{\TZVARpsi}},{\TZFy{\TZVARrhointer}{\TZVARphi}{\TZVARpsi}},{\TZFz{\TZVARrhointer}{\TZVARphi}{\TZVARpsi}}) node {$\cdot$};

\def\TZVARpsi{0}\def\TZVARphi{((pi/2)+2*pi*\TZVARh/sqrt(1+\TZVARh*\TZVARh))}
\draw[color=blue] ({\TZFx{\TZVARrhointer}{\TZVARphi}{\TZVARpsi}},{\TZFy{\TZVARrhointer}{\TZVARphi}{\TZVARpsi}},{\TZFz{\TZVARrhointer}{\TZVARphi}{\TZVARpsi}}) node {$\cdot$};

\end{scope}

\end{scope}

\begin{scope}[yshift=-2cm, xshift=3cm]

\def\TZVARrho{1}
\def\TZVARphi{pi/2}
\def\TZVARpsi{0}

\begin{scope}[param3d]

\draw[color=black!20] (-\TZVARbasesquare,-\TZVARbasesquare,0)--(\TZVARbasesquare,-\TZVARbasesquare,0)--(\TZVARbasesquare,\TZVARbasesquare,0)--(-\TZVARbasesquare,\TZVARbasesquare,0)--cycle;

\def\TZVARphi{(2*pi-0.5)}
\def\TZVARpsi{0}

\fill[color=white] plot[variable=\TZVARpsi, smooth, samples=20, domain=0:2*pi] ({\TZFx{\TZVARrho}{\TZVARphi}{\TZVARpsi}},{\TZFy{\TZVARrho}{\TZVARphi}{\TZVARpsi}},{\TZFz{\TZVARrho}{\TZVARphi}{\TZVARpsi}});

\def\TZVARpsibordone{-3.2}
\def\TZVARpsibordtwo{0.4}
\draw[color=black, fill=white, line width=0.4pt] ({\TZFx{\TZVARrho}{\TZVARphi}{\TZVARpsibordone}},{\TZFy{\TZVARrho}{\TZVARphi}{\TZVARpsibordone}},{\TZFz{\TZVARrho}{\TZVARphi}{\TZVARpsibordone}})--(0,0,0)--({\TZFx{\TZVARrho}{\TZVARphi}{\TZVARpsibordtwo}},{\TZFy{\TZVARrho}{\TZVARphi}{\TZVARpsibordtwo}},{\TZFz{\TZVARrho}{\TZVARphi}{\TZVARpsibordtwo}});

\def\TZVARpsi{0}
\def\TZVARrhointer{0.875}
\draw[color=black!40, dotted, line width=0.4pt] (0,0,0)--({\TZFx{\TZVARrhointer}{\TZVARphi}{\TZVARpsi}},{\TZFy{\TZVARrhointer}{\TZVARphi}{\TZVARpsi}},{\TZFz{\TZVARrhointer}{\TZVARphi}{\TZVARpsi}});
\draw[color=black!40, line width=0.4pt] ({\TZFx{\TZVARrhointer}{\TZVARphi}{\TZVARpsi}},{\TZFy{\TZVARrhointer}{\TZVARphi}{\TZVARpsi}},{\TZFz{\TZVARrhointer}{\TZVARphi}{\TZVARpsi}})--({\TZFx{\TZVARrho}{\TZVARphi}{\TZVARpsi}},{\TZFy{\TZVARrho}{\TZVARphi}{\TZVARpsi}},{\TZFz{\TZVARrho}{\TZVARphi}{\TZVARpsi}});

\foreach\TZVARrho in { 1}
{
\draw[color=black, line width=0.4pt, smooth, samples=20, domain=\TZVARpsibordone:\TZVARpsibordtwo] plot[variable=\TZVARpsi] ({\TZFx{\TZVARrho}{\TZVARphi}{\TZVARpsi}},{\TZFy{\TZVARrho}{\TZVARphi}{\TZVARpsi}},{\TZFz{\TZVARrho}{\TZVARphi}{\TZVARpsi}});
\draw[color=black, line width=0.4pt, smooth, samples=20, domain=\TZVARpsibordtwo:\TZVARpsibordone+2*pi] plot[variable=\TZVARpsi] ({\TZFx{\TZVARrho}{\TZVARphi}{\TZVARpsi}},{\TZFy{\TZVARrho}{\TZVARphi}{\TZVARpsi}},{\TZFz{\TZVARrho}{\TZVARphi}{\TZVARpsi}});
}

\foreach\TZVARrho in { 0.5}
{
\draw[color=red, line width=0.4pt, dotted, smooth, samples=20, domain=\TZVARpsibordone:\TZVARpsibordtwo] plot[variable=\TZVARpsi] ({\TZFx{\TZVARrho}{\TZVARphi}{\TZVARpsi}},{\TZFy{\TZVARrho}{\TZVARphi}{\TZVARpsi}},{\TZFz{\TZVARrho}{\TZVARphi}{\TZVARpsi}});
\draw[color=red, line width=0.4pt, smooth, samples=20, domain=\TZVARpsibordtwo:\TZVARpsibordone+2*pi] plot[variable=\TZVARpsi] ({\TZFx{\TZVARrho}{\TZVARphi}{\TZVARpsi}},{\TZFy{\TZVARrho}{\TZVARphi}{\TZVARpsi}},{\TZFz{\TZVARrho}{\TZVARphi}{\TZVARpsi}});
}

\def\TZVARrhointer{0.5}
\def\TZVARphiinter{4.12}
\draw[color=blue, line width=0.2pt, smooth, samples=20, domain=(pi/2):\TZVARphiinter] plot[variable=\TZVARphi] ({\TZFx{\TZVARrhointer}{\TZVARphi}{\TZVARpsi}},{\TZFy{\TZVARrhointer}{\TZVARphi}{\TZVARpsi}},{\TZFz{\TZVARrhointer}{\TZVARphi}{\TZVARpsi}});
\draw[color=blue!30, line width=0.2pt, smooth, samples=20, domain=\TZVARphiinter:\TZVARphi] plot[variable=\TZVARphi] ({\TZFx{\TZVARrhointer}{\TZVARphi}{\TZVARpsi}},{\TZFy{\TZVARrhointer}{\TZVARphi}{\TZVARpsi}},{\TZFz{\TZVARrhointer}{\TZVARphi}{\TZVARpsi}});

\def\TZVARpsi{(((pi/2)-\TZVARphi)*sqrt(1+\TZVARh*\TZVARh)/\TZVARh)}
\draw[color=red!50] ({\TZFx{\TZVARrhointer}{\TZVARphi}{\TZVARpsi}},{\TZFy{\TZVARrhointer}{\TZVARphi}{\TZVARpsi}},{\TZFz{\TZVARrhointer}{\TZVARphi}{\TZVARpsi}}) node {$\cdot$};

\def\TZVARpsi{0}\def\TZVARphi{pi/2}
\draw[color=blue] ({\TZFx{\TZVARrhointer}{\TZVARphi}{\TZVARpsi}},{\TZFy{\TZVARrhointer}{\TZVARphi}{\TZVARpsi}},{\TZFz{\TZVARrhointer}{\TZVARphi}{\TZVARpsi}}) node {$\cdot$};

\def\TZVARpsi{0}\def\TZVARphi{((pi/2)+2*pi*\TZVARh/sqrt(1+\TZVARh*\TZVARh))}
\draw[color=blue] ({\TZFx{\TZVARrhointer}{\TZVARphi}{\TZVARpsi}},{\TZFy{\TZVARrhointer}{\TZVARphi}{\TZVARpsi}},{\TZFz{\TZVARrhointer}{\TZVARphi}{\TZVARpsi}}) node {$\cdot$};
\def\TZVARpsi{0}\def\TZVARphi{((pi/2)+4*pi*\TZVARh/sqrt(1+\TZVARh*\TZVARh))}
\draw[color=blue!50] ({\TZFx{\TZVARrhointer}{\TZVARphi}{\TZVARpsi}},{\TZFy{\TZVARrhointer}{\TZVARphi}{\TZVARpsi}},{\TZFz{\TZVARrhointer}{\TZVARphi}{\TZVARpsi}}) node {$\cdot$};

\end{scope}

\end{scope}

\begin{scope}[yshift=-2cm, xshift=6cm]

\def\TZVARrho{1}
\def\TZVARphi{pi/2}
\def\TZVARpsi{0}

\begin{scope}[param3d]

\draw[color=black!20] (-\TZVARbasesquare,-\TZVARbasesquare,0)--(\TZVARbasesquare,-\TZVARbasesquare,0)--(\TZVARbasesquare,\TZVARbasesquare,0)--(-\TZVARbasesquare,\TZVARbasesquare,0)--cycle;

\def\TZVARphi{(2*pi+1.15)}
\def\TZVARpsi{0}

\fill[color=white] plot[variable=\TZVARpsi, smooth, samples=20, domain=0:2*pi] ({\TZFx{\TZVARrho}{\TZVARphi}{\TZVARpsi}},{\TZFy{\TZVARrho}{\TZVARphi}{\TZVARpsi}},{\TZFz{\TZVARrho}{\TZVARphi}{\TZVARpsi}});

\def\TZVARpsibordone{-2.75}
\def\TZVARpsibordtwo{-0.6}
\draw[color=black, fill=white, line width=0.4pt] ({\TZFx{\TZVARrho}{\TZVARphi}{\TZVARpsibordone}},{\TZFy{\TZVARrho}{\TZVARphi}{\TZVARpsibordone}},{\TZFz{\TZVARrho}{\TZVARphi}{\TZVARpsibordone}})--(0,0,0)--({\TZFx{\TZVARrho}{\TZVARphi}{\TZVARpsibordtwo}},{\TZFy{\TZVARrho}{\TZVARphi}{\TZVARpsibordtwo}},{\TZFz{\TZVARrho}{\TZVARphi}{\TZVARpsibordtwo}});

\def\TZVARpsi{0}
\def\TZVARrhointer{0.875}
\draw[color=black!40, dotted, line width=0.4pt] (0,0,0)--({\TZFx{\TZVARrhointer}{\TZVARphi}{\TZVARpsi}},{\TZFy{\TZVARrhointer}{\TZVARphi}{\TZVARpsi}},{\TZFz{\TZVARrhointer}{\TZVARphi}{\TZVARpsi}});
\draw[color=black!40, line width=0.4pt] ({\TZFx{\TZVARrhointer}{\TZVARphi}{\TZVARpsi}},{\TZFy{\TZVARrhointer}{\TZVARphi}{\TZVARpsi}},{\TZFz{\TZVARrhointer}{\TZVARphi}{\TZVARpsi}})--({\TZFx{\TZVARrho}{\TZVARphi}{\TZVARpsi}},{\TZFy{\TZVARrho}{\TZVARphi}{\TZVARpsi}},{\TZFz{\TZVARrho}{\TZVARphi}{\TZVARpsi}});

\foreach\TZVARrho in { 1}
{
\draw[color=black, line width=0.4pt, smooth, samples=20, domain=\TZVARpsibordone:\TZVARpsibordtwo] plot[variable=\TZVARpsi] ({\TZFx{\TZVARrho}{\TZVARphi}{\TZVARpsi}},{\TZFy{\TZVARrho}{\TZVARphi}{\TZVARpsi}},{\TZFz{\TZVARrho}{\TZVARphi}{\TZVARpsi}});
\draw[color=black, line width=0.4pt, smooth, samples=20, domain=\TZVARpsibordtwo:\TZVARpsibordone+2*pi] plot[variable=\TZVARpsi] ({\TZFx{\TZVARrho}{\TZVARphi}{\TZVARpsi}},{\TZFy{\TZVARrho}{\TZVARphi}{\TZVARpsi}},{\TZFz{\TZVARrho}{\TZVARphi}{\TZVARpsi}});
}

\foreach\TZVARrho in { 0.5}
{
\draw[color=red, line width=0.4pt, smooth, samples=20, domain=\TZVARpsibordone:\TZVARpsibordtwo] plot[variable=\TZVARpsi] ({\TZFx{\TZVARrho}{\TZVARphi}{\TZVARpsi}},{\TZFy{\TZVARrho}{\TZVARphi}{\TZVARpsi}},{\TZFz{\TZVARrho}{\TZVARphi}{\TZVARpsi}});
\draw[color=red, line width=0.4pt, dotted, smooth, samples=20, domain=\TZVARpsibordtwo:\TZVARpsibordone+2*pi] plot[variable=\TZVARpsi] ({\TZFx{\TZVARrho}{\TZVARphi}{\TZVARpsi}},{\TZFy{\TZVARrho}{\TZVARphi}{\TZVARpsi}},{\TZFz{\TZVARrho}{\TZVARphi}{\TZVARpsi}});
}

\def\TZVARrhointer{0.5}
\def\TZVARphiinter{3.3}
\draw[color=blue!30, line width=0.3pt, smooth, samples=20, domain=(pi/2):\TZVARphiinter] plot[variable=\TZVARphi] ({\TZFx{\TZVARrhointer}{\TZVARphi}{\TZVARpsi}},{\TZFy{\TZVARrhointer}{\TZVARphi}{\TZVARpsi}},{\TZFz{\TZVARrhointer}{\TZVARphi}{\TZVARpsi}});
\draw[color=blue, line width=0.3pt, smooth, samples=20, domain=\TZVARphiinter:\TZVARphi] plot[variable=\TZVARphi] ({\TZFx{\TZVARrhointer}{\TZVARphi}{\TZVARpsi}},{\TZFy{\TZVARrhointer}{\TZVARphi}{\TZVARpsi}},{\TZFz{\TZVARrhointer}{\TZVARphi}{\TZVARpsi}});

\def\TZVARpsi{(((pi/2)-\TZVARphi)*sqrt(1+\TZVARh*\TZVARh)/\TZVARh)}
\draw[color=red] ({\TZFx{\TZVARrhointer}{\TZVARphi}{\TZVARpsi}},{\TZFy{\TZVARrhointer}{\TZVARphi}{\TZVARpsi}},{\TZFz{\TZVARrhointer}{\TZVARphi}{\TZVARpsi}}) node {$\cdot$};

\def\TZVARpsi{0}\def\TZVARphi{pi/2}
\draw[color=blue!50] ({\TZFx{\TZVARrhointer}{\TZVARphi}{\TZVARpsi}},{\TZFy{\TZVARrhointer}{\TZVARphi}{\TZVARpsi}},{\TZFz{\TZVARrhointer}{\TZVARphi}{\TZVARpsi}}) node {$\cdot$};

\def\TZVARpsi{0}\def\TZVARphi{((pi/2)+2*pi*\TZVARh/sqrt(1+\TZVARh*\TZVARh))}
\draw[color=blue] ({\TZFx{\TZVARrhointer}{\TZVARphi}{\TZVARpsi}},{\TZFy{\TZVARrhointer}{\TZVARphi}{\TZVARpsi}},{\TZFz{\TZVARrhointer}{\TZVARphi}{\TZVARpsi}}) node {$\cdot$};
\def\TZVARpsi{0}\def\TZVARphi{((pi/2)+4*pi*\TZVARh/sqrt(1+\TZVARh*\TZVARh))}
\draw[color=blue] ({\TZFx{\TZVARrhointer}{\TZVARphi}{\TZVARpsi}},{\TZFy{\TZVARrhointer}{\TZVARphi}{\TZVARpsi}},{\TZFz{\TZVARrhointer}{\TZVARphi}{\TZVARpsi}}) node {$\cdot$};
\def\TZVARpsi{0}\def\TZVARphi{((pi/2)+6*pi*\TZVARh/sqrt(1+\TZVARh*\TZVARh))}
\draw[color=blue] ({\TZFx{\TZVARrhointer}{\TZVARphi}{\TZVARpsi}},{\TZFy{\TZVARrhointer}{\TZVARphi}{\TZVARpsi}},{\TZFz{\TZVARrhointer}{\TZVARphi}{\TZVARpsi}}) node {$\cdot$};

\end{scope}

\end{scope}

\end{tikzpicture}
\end{figure}
}%

When we want to transcribe this game in mathematical terms, we change our point of view; instead of being a spectator that sees a ball rolling on the ground we become someone living on the (surface of the) ball that sees the ground (the tangent space) turning around him. 
Said differently, we no longer see the ball rolling on the ground, we see the ground (\ie, the tangent space) rolling on the ball. 
An ant on the ball who thinks the ball is the center of the universe and thus (necessarily) fixed, sees the ground moving around him. 
This is the movement of the rolling tangent space. 
And of course, we replace the ball by an arbitrary submanifold $M$ of $\RR^\nu$. 
\WelNietTikZpicture
{%
\begin{figure}[!ht]
\begin{tikzpicture}[scale=1.0, line cap=round,line join=round,x=\TZschaal, y=\TZschaal]

\def\TZemme{0.4}
\def\TZhoek{0}

\def\TZVARsep{2.75cm}

\def\TZVARr{0.5} 
\def\TZVARrinter{0.25}
\def\TZVARh{0.3}
\def\TZVARalpha{(\TZVARh/sqrt(1+(\TZVARh)*(\TZVARh)))} 

\def\TZVARptONEx{0.6*\TZVARr}
\def\TZVARptONEy{(0.6*\TZVARr/\TZVARalpha)}
\def\TZVARptONEz{(-\TZVARptONEx/\TZVARh)}

\let\TZVARptTWOx\TZVARptONEx
\def\TZVARptTWOy{(-0.1*\TZVARr/\TZVARalpha)}
\def\TZVARptTWOz{(-\TZVARptTWOx/\TZVARh)}

\def\TZVARptTHREEx{-0.2*\TZVARr}
\let\TZVARptTHREEy\TZVARptTWOy
\def\TZVARptTHREEz{(-\TZVARptTHREEx/\TZVARh)}

\let\TZVARptFOURx\TZVARptTHREEx
\let\TZVARptFOURy\TZVARptONEy
\def\TZVARptFOURz{(-\TZVARptFOURx/\TZVARh)}

\def\TZColorCone{gray!30}

\def\TZFBx#1#2#3{%
(
(1+\TZVARh*\TZVARh*cos(\TZVARalpha*(\TZVARt-\TZVARs) r))*(#1)/(1+\TZVARh*\TZVARh)
+
\TZVARalpha*sin(\TZVARalpha*(\TZVARt-\TZVARs) r)*(#2)
+
\TZVARh*(1-cos(\TZVARalpha*(\TZVARt-\TZVARs) r))*(#3)/(1+\TZVARh*\TZVARh)
)}
\def\TZFBy#1#2#3{%
(
-\TZVARalpha*sin(\TZVARalpha*(\TZVARt-\TZVARs) r)*(#1)
+
cos(\TZVARalpha*(\TZVARt-\TZVARs) r)*(#2)
+
sin(\TZVARalpha*(\TZVARt-\TZVARs) r)*(#3)/sqrt(1+\TZVARh*\TZVARh)
)}
\def\TZFBz#1#2#3{%
(
\TZVARh*(1-cos(\TZVARalpha*(\TZVARt-\TZVARs) r))*(#1)/(1+\TZVARh*\TZVARh)
-
sin(\TZVARalpha*(\TZVARt-\TZVARs) r)*(#2)/sqrt(1+\TZVARh*\TZVARh)
+
(\TZVARh*\TZVARh+cos(\TZVARalpha*(\TZVARt-\TZVARs) r))*(#3)/(1+\TZVARh*\TZVARh)
)}

\def\TZFAx#1#2#3{%
(
cos(\TZVARt r)*\TZFBx{#1}{#2}{#3} - sin(\TZVARt r)*\TZFBy{#1}{#2}{#3}
)}
\def\TZFAy#1#2#3{%
(
sin(\TZVARt r)*\TZFBx{#1}{#2}{#3} + cos(\TZVARt r)*\TZFBy{#1}{#2}{#3}
)}
\let\TZFAz\TZFBz

\begin{scope}[xshift=1cm, yshift=-3cm]

\begin{scope}[yshift=0.75cm, xshift=-0cm, param3d]
\draw[color=black, line width=0.15pt] plot[variable=\TZVARt, smooth, samples=20, domain=-pi/2:pi/2] ({\TZVARrinter*cos(\TZVARt r)},{\TZVARrinter*sin(\TZVARt r)},{-\TZVARrinter/\TZVARh});
\draw[color=black, dotted, line width=0.15pt, smooth, samples=20, domain=pi/2:3*pi/2] plot[variable=\TZVARt] ({\TZVARrinter*cos(\TZVARt r)},{\TZVARrinter*sin(\TZVARt r)},{-\TZVARrinter/\TZVARh});
\end{scope}
\draw (0.5,0) node[right] {\footnotesize the initial curve (a circle) on which the tangent space rolls};

\begin{scope}[yshift=-0.5cm]
\draw[color=red, line width=0.4pt] (-0.5,0)--(0.5,0);
\draw (0.5,0) node[right] {\footnotesize its trace on the tangent space};
\end{scope}

\begin{scope}[yshift=-1cm]
\draw[color=red, dashed, line width=0.4pt] (-0.5,0)--(0.5,0);
\draw (0.5,0) node[right] {\footnotesize the contact line of the cone with the rolling tangent space};
\end{scope}

\end{scope} 

\begin{scope}[param3d]

\draw[color=black, line width=0.4pt, fill=\TZColorCone] plot[variable=\TZVARt, smooth, samples=20, domain=-pi/2:pi/2] ({\TZVARr*cos(\TZVARt r)},{\TZVARr*sin(\TZVARt r)},{-\TZVARr/\TZVARh})--(0,0,0)--cycle;
\draw[color=black, dotted, line width=0.4pt, smooth, samples=20, domain=pi/2:3*pi/2] plot[variable=\TZVARt] ({\TZVARr*cos(\TZVARt r)},{\TZVARr*sin(\TZVARt r)},{-\TZVARr/\TZVARh});
\draw[color=black, line width=0.15pt] plot[variable=\TZVARt, smooth, samples=20, domain=-pi/2:pi/2] ({\TZVARrinter*cos(\TZVARt r)},{\TZVARrinter*sin(\TZVARt r)},{-\TZVARrinter/\TZVARh});
\draw[color=black, dotted, line width=0.15pt, smooth, samples=20, domain=pi/2:3*pi/2] plot[variable=\TZVARt] ({\TZVARrinter*cos(\TZVARt r)},{\TZVARrinter*sin(\TZVARt r)},{-\TZVARrinter/\TZVARh});

\def\TZVARs{0}

\foreach \TZVARt in {0}
{
\draw[color=brown, line width=0.4pt, fill=white, opacity=0.9] 
({\TZFAx{\TZVARptONEx}{\TZVARptONEy}{\TZVARptONEz}},
{\TZFAy{\TZVARptONEx}{\TZVARptONEy}{\TZVARptONEz}}, 
{\TZFAz{\TZVARptONEx}{\TZVARptONEy}{\TZVARptONEz}})
--
({\TZFAx{\TZVARptTWOx}{\TZVARptTWOy}{\TZVARptTWOz}},
{\TZFAy{\TZVARptTWOx}{\TZVARptTWOy}{\TZVARptTWOz}}, 
{\TZFAz{\TZVARptTWOx}{\TZVARptTWOy}{\TZVARptTWOz}})
--
({\TZFAx{\TZVARptTHREEx}{\TZVARptTHREEy}{\TZVARptTHREEz}},
{\TZFAy{\TZVARptTHREEx}{\TZVARptTHREEy}{\TZVARptTHREEz}}, 
{\TZFAz{\TZVARptTHREEx}{\TZVARptTHREEy}{\TZVARptTHREEz}})
--
({\TZFAx{\TZVARptFOURx}{\TZVARptFOURy}{\TZVARptFOURz}},
{\TZFAy{\TZVARptFOURx}{\TZVARptFOURy}{\TZVARptFOURz}}, 
{\TZFAz{\TZVARptFOURx}{\TZVARptFOURy}{\TZVARptFOURz}})
--cycle;
\draw[color=red, line width=0.4pt, smooth, samples=20, domain=0:\TZVARt] plot[variable=\TZVARs] ({\TZFAx{\TZVARrinter}{0}{-\TZVARrinter/\TZVARh}}, {\TZFAy{\TZVARrinter}{0}{-\TZVARrinter/\TZVARh}}, {\TZFAz{\TZVARrinter}{0}{-\TZVARrinter/\TZVARh}});
\def\TZVARtouchlen{0.6}
\draw[color=red, dashed, line width=0.2pt] (0,0,0)--({\TZVARtouchlen*\TZVARr*cos(\TZVARt)},{-\TZVARtouchlen*\TZVARr*sin(\TZVARt)}, {-\TZVARtouchlen*\TZVARr/\TZVARh});
}

\end{scope}
\draw (0,-1.85) node[below] {\scriptsize $\varphi=0$};

\begin{scope}[xshift=\TZVARsep]

\begin{scope}[param3d]

\draw[color=black, line width=0.4pt, fill=\TZColorCone] plot[variable=\TZVARt, smooth, samples=20, domain=-pi/2:pi/2] ({\TZVARr*cos(\TZVARt r)},{\TZVARr*sin(\TZVARt r)},{-\TZVARr/\TZVARh})--(0,0,0)--cycle;
\draw[color=black, dotted, line width=0.4pt, smooth, samples=20, domain=pi/2:3*pi/2] plot[variable=\TZVARt] ({\TZVARr*cos(\TZVARt r)},{\TZVARr*sin(\TZVARt r)},{-\TZVARr/\TZVARh});

\draw[color=black, line width=0.15pt] plot[variable=\TZVARt, smooth, samples=20, domain=-pi/2:pi/2] ({\TZVARrinter*cos(\TZVARt r)},{\TZVARrinter*sin(\TZVARt r)},{-\TZVARrinter/\TZVARh});
\draw[color=black, dotted, line width=0.15pt, smooth, samples=20, domain=pi/2:3*pi/2] plot[variable=\TZVARt] ({\TZVARrinter*cos(\TZVARt r)},{\TZVARrinter*sin(\TZVARt r)},{-\TZVARrinter/\TZVARh});

\def\TZVARs{0}

\foreach \TZVARt in {0.8}
{
\draw[color=brown, line width=0.4pt, fill=white, opacity=0.9] 
({\TZFAx{\TZVARptONEx}{\TZVARptONEy}{\TZVARptONEz}},
{\TZFAy{\TZVARptONEx}{\TZVARptONEy}{\TZVARptONEz}}, 
{\TZFAz{\TZVARptONEx}{\TZVARptONEy}{\TZVARptONEz}})
--
({\TZFAx{\TZVARptTWOx}{\TZVARptTWOy}{\TZVARptTWOz}},
{\TZFAy{\TZVARptTWOx}{\TZVARptTWOy}{\TZVARptTWOz}}, 
{\TZFAz{\TZVARptTWOx}{\TZVARptTWOy}{\TZVARptTWOz}})
--
({\TZFAx{\TZVARptTHREEx}{\TZVARptTHREEy}{\TZVARptTHREEz}},
{\TZFAy{\TZVARptTHREEx}{\TZVARptTHREEy}{\TZVARptTHREEz}}, 
{\TZFAz{\TZVARptTHREEx}{\TZVARptTHREEy}{\TZVARptTHREEz}})
--
({\TZFAx{\TZVARptFOURx}{\TZVARptFOURy}{\TZVARptFOURz}},
{\TZFAy{\TZVARptFOURx}{\TZVARptFOURy}{\TZVARptFOURz}}, 
{\TZFAz{\TZVARptFOURx}{\TZVARptFOURy}{\TZVARptFOURz}})
--cycle;
\draw[color=red, line width=0.4pt, smooth, samples=20, domain=0:\TZVARt] plot[variable=\TZVARs] ({\TZFAx{\TZVARrinter}{0}{-\TZVARrinter/\TZVARh}}, {\TZFAy{\TZVARrinter}{0}{-\TZVARrinter/\TZVARh}}, {\TZFAz{\TZVARrinter}{0}{-\TZVARrinter/\TZVARh}});
\def\TZVARtouchlen{0.61}
\draw[color=red, dashed, line width=0.2pt] (0,0,0)--({\TZVARtouchlen*\TZVARr*cos(\TZVARt r)},{\TZVARtouchlen*\TZVARr*sin(\TZVARt r)}, {-\TZVARtouchlen*\TZVARr/\TZVARh});
}

\end{scope}
\draw (0,-1.85) node[below] {\scriptsize $\varphi=0.8$};
\end{scope}

\begin{scope}[xshift=2*\TZVARsep]

\begin{scope}[param3d]

\draw[color=black, line width=0.4pt, fill=\TZColorCone] plot[variable=\TZVARt, smooth, samples=20, domain=-pi/2:pi/2] ({\TZVARr*cos(\TZVARt r)},{\TZVARr*sin(\TZVARt r)},{-\TZVARr/\TZVARh})--(0,0,0)--cycle;
\draw[color=black, dotted, line width=0.4pt, smooth, samples=20, domain=pi/2:3*pi/2] plot[variable=\TZVARt] ({\TZVARr*cos(\TZVARt r)},{\TZVARr*sin(\TZVARt r)},{-\TZVARr/\TZVARh});

\draw[color=black, line width=0.15pt] plot[variable=\TZVARt, smooth, samples=20, domain=-pi/2:pi/2] ({\TZVARrinter*cos(\TZVARt r)},{\TZVARrinter*sin(\TZVARt r)},{-\TZVARrinter/\TZVARh});
\draw[color=black, dotted, line width=0.15pt, smooth, samples=20, domain=pi/2:3*pi/2] plot[variable=\TZVARt] ({\TZVARrinter*cos(\TZVARt r)},{\TZVARrinter*sin(\TZVARt r)},{-\TZVARrinter/\TZVARh});

\def\TZVARs{0}

\foreach \TZVARt in {(1.3)}
{
\draw[color=brown, line width=0.4pt, fill=white, opacity=0.8] 
({\TZFAx{\TZVARptONEx}{\TZVARptONEy}{\TZVARptONEz}},
{\TZFAy{\TZVARptONEx}{\TZVARptONEy}{\TZVARptONEz}}, 
{\TZFAz{\TZVARptONEx}{\TZVARptONEy}{\TZVARptONEz}})
--
({\TZFAx{\TZVARptTWOx}{\TZVARptTWOy}{\TZVARptTWOz}},
{\TZFAy{\TZVARptTWOx}{\TZVARptTWOy}{\TZVARptTWOz}}, 
{\TZFAz{\TZVARptTWOx}{\TZVARptTWOy}{\TZVARptTWOz}})
--
({\TZFAx{\TZVARptTHREEx}{\TZVARptTHREEy}{\TZVARptTHREEz}},
{\TZFAy{\TZVARptTHREEx}{\TZVARptTHREEy}{\TZVARptTHREEz}}, 
{\TZFAz{\TZVARptTHREEx}{\TZVARptTHREEy}{\TZVARptTHREEz}})
--
({\TZFAx{\TZVARptFOURx}{\TZVARptFOURy}{\TZVARptFOURz}},
{\TZFAy{\TZVARptFOURx}{\TZVARptFOURy}{\TZVARptFOURz}}, 
{\TZFAz{\TZVARptFOURx}{\TZVARptFOURy}{\TZVARptFOURz}})
--cycle;
\draw[color=red, line width=0.4pt, smooth, samples=20, domain=0:\TZVARt] plot[variable=\TZVARs] ({\TZFAx{\TZVARrinter}{0}{-\TZVARrinter/\TZVARh}}, {\TZFAy{\TZVARrinter}{0}{-\TZVARrinter/\TZVARh}}, {\TZFAz{\TZVARrinter}{0}{-\TZVARrinter/\TZVARh}});
\def\TZVARtouchlen{0.64}
\draw[color=red, dashed, line width=0.2pt] (0,0,0)--({\TZVARtouchlen*\TZVARr*cos(\TZVARt r)},{\TZVARtouchlen*\TZVARr*sin(\TZVARt r)}, {-\TZVARtouchlen*\TZVARr/\TZVARh});
}

\end{scope}
\draw (0,-1.85) node[below] {\scriptsize $\varphi=1.3$};
\end{scope}

\begin{scope}[xshift=3*\TZVARsep]

\begin{scope}[param3d]

\def\TZVARs{0}

\foreach \TZVARt in {(2.1)}
{
\draw[color=brown, line width=0.4pt, fill=white] 
({\TZFAx{\TZVARptONEx}{\TZVARptONEy}{\TZVARptONEz}},
{\TZFAy{\TZVARptONEx}{\TZVARptONEy}{\TZVARptONEz}}, 
{\TZFAz{\TZVARptONEx}{\TZVARptONEy}{\TZVARptONEz}})
--
({\TZFAx{\TZVARptTWOx}{\TZVARptTWOy}{\TZVARptTWOz}},
{\TZFAy{\TZVARptTWOx}{\TZVARptTWOy}{\TZVARptTWOz}}, 
{\TZFAz{\TZVARptTWOx}{\TZVARptTWOy}{\TZVARptTWOz}})
--
({\TZFAx{\TZVARptTHREEx}{\TZVARptTHREEy}{\TZVARptTHREEz}},
{\TZFAy{\TZVARptTHREEx}{\TZVARptTHREEy}{\TZVARptTHREEz}}, 
{\TZFAz{\TZVARptTHREEx}{\TZVARptTHREEy}{\TZVARptTHREEz}})
--
({\TZFAx{\TZVARptFOURx}{\TZVARptFOURy}{\TZVARptFOURz}},
{\TZFAy{\TZVARptFOURx}{\TZVARptFOURy}{\TZVARptFOURz}}, 
{\TZFAz{\TZVARptFOURx}{\TZVARptFOURy}{\TZVARptFOURz}})
--cycle;
\draw[color=red, line width=0.4pt, smooth, samples=20, domain=0:\TZVARt] plot[variable=\TZVARs] ({\TZFAx{\TZVARrinter}{0}{-\TZVARrinter/\TZVARh}}, {\TZFAy{\TZVARrinter}{0}{-\TZVARrinter/\TZVARh}}, {\TZFAz{\TZVARrinter}{0}{-\TZVARrinter/\TZVARh}});
\def\TZVARtouchlen{0.73} 
\draw[color=red, dashed, line width=0.2pt] (0,0,0)--({\TZVARtouchlen*\TZVARr*cos(\TZVARt r)},{\TZVARtouchlen*\TZVARr*sin(\TZVARt r)}, {-\TZVARtouchlen*\TZVARr/\TZVARh});
}

\draw[color=black, line width=0.4pt, fill=\TZColorCone, opacity=0.8] plot[variable=\TZVARt, smooth, samples=20, domain=-pi/2:pi/2] ({\TZVARr*cos(\TZVARt r)},{\TZVARr*sin(\TZVARt r)},{-\TZVARr/\TZVARh})--(0,0,0)--cycle;
\draw[color=black, dotted, line width=0.4pt, smooth, samples=20, domain=pi/2:3*pi/2] plot[variable=\TZVARt] ({\TZVARr*cos(\TZVARt r)},{\TZVARr*sin(\TZVARt r)},{-\TZVARr/\TZVARh});

\draw[color=black, line width=0.15pt] plot[variable=\TZVARt, smooth, samples=20, domain=-pi/2:pi/2] ({\TZVARrinter*cos(\TZVARt r)},{\TZVARrinter*sin(\TZVARt r)},{-\TZVARrinter/\TZVARh});
\draw[color=black, dotted, line width=0.15pt, smooth, samples=20, domain=pi/2:3*pi/2] plot[variable=\TZVARt] ({\TZVARrinter*cos(\TZVARt r)},{\TZVARrinter*sin(\TZVARt r)},{-\TZVARrinter/\TZVARh});

\end{scope}
\draw (0,-1.85) node[below] {\scriptsize $\varphi=2.1$};
\end{scope}

\begin{scope}[xshift=4*\TZVARsep]

\begin{scope}[param3d]

\def\TZVARs{0}

\foreach \TZVARt in {(2.9)}
{
\draw[color=brown, line width=0.4pt, fill=white] 
({\TZFAx{\TZVARptONEx}{\TZVARptONEy}{\TZVARptONEz}},
{\TZFAy{\TZVARptONEx}{\TZVARptONEy}{\TZVARptONEz}}, 
{\TZFAz{\TZVARptONEx}{\TZVARptONEy}{\TZVARptONEz}})
--
({\TZFAx{\TZVARptTWOx}{\TZVARptTWOy}{\TZVARptTWOz}},
{\TZFAy{\TZVARptTWOx}{\TZVARptTWOy}{\TZVARptTWOz}}, 
{\TZFAz{\TZVARptTWOx}{\TZVARptTWOy}{\TZVARptTWOz}})
--
({\TZFAx{\TZVARptTHREEx}{\TZVARptTHREEy}{\TZVARptTHREEz}},
{\TZFAy{\TZVARptTHREEx}{\TZVARptTHREEy}{\TZVARptTHREEz}}, 
{\TZFAz{\TZVARptTHREEx}{\TZVARptTHREEy}{\TZVARptTHREEz}})
--
({\TZFAx{\TZVARptFOURx}{\TZVARptFOURy}{\TZVARptFOURz}},
{\TZFAy{\TZVARptFOURx}{\TZVARptFOURy}{\TZVARptFOURz}}, 
{\TZFAz{\TZVARptFOURx}{\TZVARptFOURy}{\TZVARptFOURz}})
--cycle;
\draw[color=red, line width=0.4pt, smooth, samples=20, domain=0:\TZVARt] plot[variable=\TZVARs] ({\TZFAx{\TZVARrinter}{0}{-\TZVARrinter/\TZVARh}}, {\TZFAy{\TZVARrinter}{0}{-\TZVARrinter/\TZVARh}}, {\TZFAz{\TZVARrinter}{0}{-\TZVARrinter/\TZVARh}});
\def\TZVARtouchlen{0.805}
\draw[color=red, dashed, line width=0.2pt] (0,0,0)--({\TZVARtouchlen*\TZVARr*cos(\TZVARt r)},{\TZVARtouchlen*\TZVARr*sin(\TZVARt r)}, {-\TZVARtouchlen*\TZVARr/\TZVARh});
}

\draw[color=black, line width=0.4pt, fill=\TZColorCone, opacity=0.8] plot[variable=\TZVARt, smooth, samples=20, domain=-pi/2:pi/2] ({\TZVARr*cos(\TZVARt r)},{\TZVARr*sin(\TZVARt r)},{-\TZVARr/\TZVARh})--(0,0,0)--cycle;
\draw[color=black, dotted, line width=0.4pt, smooth, samples=20, domain=pi/2:3*pi/2] plot[variable=\TZVARt] ({\TZVARr*cos(\TZVARt r)},{\TZVARr*sin(\TZVARt r)},{-\TZVARr/\TZVARh});

\draw[color=black, line width=0.15pt] plot[variable=\TZVARt, smooth, samples=20, domain=-pi/2:pi/2] ({\TZVARrinter*cos(\TZVARt r)},{\TZVARrinter*sin(\TZVARt r)},{-\TZVARrinter/\TZVARh});
\draw[color=black, dotted, line width=0.15pt, smooth, samples=20, domain=pi/2:3*pi/2] plot[variable=\TZVARt] ({\TZVARrinter*cos(\TZVARt r)},{\TZVARrinter*sin(\TZVARt r)},{-\TZVARrinter/\TZVARh});

\end{scope}
\draw (0,-1.85) node[below] {\scriptsize $\varphi=2.9$};
\end{scope}

\begin{scope}[yshift=-5cm]

\begin{scope}

\begin{scope}[param3d]

\def\TZVARs{0}

\foreach \TZVARt in {(3.6)}
{
\draw[color=brown, line width=0.4pt, fill=white] 
({\TZFAx{\TZVARptONEx}{\TZVARptONEy}{\TZVARptONEz}},
{\TZFAy{\TZVARptONEx}{\TZVARptONEy}{\TZVARptONEz}}, 
{\TZFAz{\TZVARptONEx}{\TZVARptONEy}{\TZVARptONEz}})
--
({\TZFAx{\TZVARptTWOx}{\TZVARptTWOy}{\TZVARptTWOz}},
{\TZFAy{\TZVARptTWOx}{\TZVARptTWOy}{\TZVARptTWOz}}, 
{\TZFAz{\TZVARptTWOx}{\TZVARptTWOy}{\TZVARptTWOz}})
--
({\TZFAx{\TZVARptTHREEx}{\TZVARptTHREEy}{\TZVARptTHREEz}},
{\TZFAy{\TZVARptTHREEx}{\TZVARptTHREEy}{\TZVARptTHREEz}}, 
{\TZFAz{\TZVARptTHREEx}{\TZVARptTHREEy}{\TZVARptTHREEz}})
--
({\TZFAx{\TZVARptFOURx}{\TZVARptFOURy}{\TZVARptFOURz}},
{\TZFAy{\TZVARptFOURx}{\TZVARptFOURy}{\TZVARptFOURz}}, 
{\TZFAz{\TZVARptFOURx}{\TZVARptFOURy}{\TZVARptFOURz}})
--cycle;
\draw[color=red, line width=0.4pt, smooth, samples=20, domain=0:\TZVARt] plot[variable=\TZVARs] ({\TZFAx{\TZVARrinter}{0}{-\TZVARrinter/\TZVARh}}, {\TZFAy{\TZVARrinter}{0}{-\TZVARrinter/\TZVARh}}, {\TZFAz{\TZVARrinter}{0}{-\TZVARrinter/\TZVARh}});
\def\TZVARtouchlen{0.69}
\draw[color=red, dashed, line width=0.2pt] (0,0,0)--({\TZVARtouchlen*\TZVARr*cos(\TZVARt r)},{\TZVARtouchlen*\TZVARr*sin(\TZVARt r)}, {-\TZVARtouchlen*\TZVARr/\TZVARh});
}

\draw[color=black, line width=0.4pt, fill=\TZColorCone, opacity=0.8] plot[variable=\TZVARt, smooth, samples=20, domain=-pi/2:pi/2] ({\TZVARr*cos(\TZVARt r)},{\TZVARr*sin(\TZVARt r)},{-\TZVARr/\TZVARh})--(0,0,0)--cycle;
\draw[color=black, dotted, line width=0.4pt, smooth, samples=20, domain=pi/2:3*pi/2] plot[variable=\TZVARt] ({\TZVARr*cos(\TZVARt r)},{\TZVARr*sin(\TZVARt r)},{-\TZVARr/\TZVARh});

\draw[color=black, line width=0.15pt] plot[variable=\TZVARt, smooth, samples=20, domain=-pi/2:pi/2] ({\TZVARrinter*cos(\TZVARt r)},{\TZVARrinter*sin(\TZVARt r)},{-\TZVARrinter/\TZVARh});
\draw[color=black, dotted, line width=0.15pt, smooth, samples=20, domain=pi/2:3*pi/2] plot[variable=\TZVARt] ({\TZVARrinter*cos(\TZVARt r)},{\TZVARrinter*sin(\TZVARt r)},{-\TZVARrinter/\TZVARh});

\end{scope}
\draw (0,-1.85) node[below] {\scriptsize $\varphi=3.6$};
\end{scope}

\begin{scope}[xshift=\TZVARsep]

\begin{scope}[param3d]

\def\TZVARs{0}

\foreach \TZVARt in {(4.3)}
{
\draw[color=brown, line width=0.4pt, fill=white] 
({\TZFAx{\TZVARptONEx}{\TZVARptONEy}{\TZVARptONEz}},
{\TZFAy{\TZVARptONEx}{\TZVARptONEy}{\TZVARptONEz}}, 
{\TZFAz{\TZVARptONEx}{\TZVARptONEy}{\TZVARptONEz}})
--
({\TZFAx{\TZVARptTWOx}{\TZVARptTWOy}{\TZVARptTWOz}},
{\TZFAy{\TZVARptTWOx}{\TZVARptTWOy}{\TZVARptTWOz}}, 
{\TZFAz{\TZVARptTWOx}{\TZVARptTWOy}{\TZVARptTWOz}})
--
({\TZFAx{\TZVARptTHREEx}{\TZVARptTHREEy}{\TZVARptTHREEz}},
{\TZFAy{\TZVARptTHREEx}{\TZVARptTHREEy}{\TZVARptTHREEz}}, 
{\TZFAz{\TZVARptTHREEx}{\TZVARptTHREEy}{\TZVARptTHREEz}})
--
({\TZFAx{\TZVARptFOURx}{\TZVARptFOURy}{\TZVARptFOURz}},
{\TZFAy{\TZVARptFOURx}{\TZVARptFOURy}{\TZVARptFOURz}}, 
{\TZFAz{\TZVARptFOURx}{\TZVARptFOURy}{\TZVARptFOURz}})
--cycle;
\draw[color=red, line width=0.4pt, smooth, samples=20, domain=0:\TZVARt] plot[variable=\TZVARs] ({\TZFAx{\TZVARrinter}{0}{-\TZVARrinter/\TZVARh}}, {\TZFAy{\TZVARrinter}{0}{-\TZVARrinter/\TZVARh}}, {\TZFAz{\TZVARrinter}{0}{-\TZVARrinter/\TZVARh}});
\def\TZVARtouchlen{0.63}
\draw[color=red, dashed, line width=0.2pt] (0,0,0)--({\TZVARtouchlen*\TZVARr*cos(\TZVARt r)},{\TZVARtouchlen*\TZVARr*sin(\TZVARt r)}, {-\TZVARtouchlen*\TZVARr/\TZVARh});
}

\draw[color=black, line width=0.4pt, fill=\TZColorCone, opacity=0.8] plot[variable=\TZVARt, smooth, samples=20, domain=-pi/2:pi/2] ({\TZVARr*cos(\TZVARt r)},{\TZVARr*sin(\TZVARt r)},{-\TZVARr/\TZVARh})--(0,0,0)--cycle;
\draw[color=black, dotted, line width=0.4pt, smooth, samples=20, domain=pi/2:3*pi/2] plot[variable=\TZVARt] ({\TZVARr*cos(\TZVARt r)},{\TZVARr*sin(\TZVARt r)},{-\TZVARr/\TZVARh});

\draw[color=black, line width=0.15pt] plot[variable=\TZVARt, smooth, samples=20, domain=-pi/2:pi/2] ({\TZVARrinter*cos(\TZVARt r)},{\TZVARrinter*sin(\TZVARt r)},{-\TZVARrinter/\TZVARh});
\draw[color=black, dotted, line width=0.15pt, smooth, samples=20, domain=pi/2:3*pi/2] plot[variable=\TZVARt] ({\TZVARrinter*cos(\TZVARt r)},{\TZVARrinter*sin(\TZVARt r)},{-\TZVARrinter/\TZVARh});

\end{scope}
\draw (0,-1.85) node[below] {\scriptsize $\varphi=4.3$};
\end{scope}

\begin{scope}[xshift=2*\TZVARsep]

\begin{scope}[param3d]

\draw[color=black, line width=0.4pt, fill=\TZColorCone] plot[variable=\TZVARt, smooth, samples=20, domain=-pi/2:pi/2] ({\TZVARr*cos(\TZVARt r)},{\TZVARr*sin(\TZVARt r)},{-\TZVARr/\TZVARh})--(0,0,0)--cycle;
\draw[color=black, dotted, line width=0.4pt, smooth, samples=20, domain=pi/2:3*pi/2] plot[variable=\TZVARt] ({\TZVARr*cos(\TZVARt r)},{\TZVARr*sin(\TZVARt r)},{-\TZVARr/\TZVARh});

\draw[color=black, line width=0.15pt] plot[variable=\TZVARt, smooth, samples=20, domain=-pi/2:pi/2] ({\TZVARrinter*cos(\TZVARt r)},{\TZVARrinter*sin(\TZVARt r)},{-\TZVARrinter/\TZVARh});
\draw[color=black, dotted, line width=0.15pt, smooth, samples=20, domain=pi/2:3*pi/2] plot[variable=\TZVARt] ({\TZVARrinter*cos(\TZVARt r)},{\TZVARrinter*sin(\TZVARt r)},{-\TZVARrinter/\TZVARh});

\def\TZVARs{0}

\foreach \TZVARt in {(4.8)}
{
\draw[color=brown, line width=0.4pt, fill=white, opacity=0.8] 
({\TZFAx{\TZVARptONEx}{\TZVARptONEy}{\TZVARptONEz}},
{\TZFAy{\TZVARptONEx}{\TZVARptONEy}{\TZVARptONEz}}, 
{\TZFAz{\TZVARptONEx}{\TZVARptONEy}{\TZVARptONEz}})
--
({\TZFAx{\TZVARptTWOx}{\TZVARptTWOy}{\TZVARptTWOz}},
{\TZFAy{\TZVARptTWOx}{\TZVARptTWOy}{\TZVARptTWOz}}, 
{\TZFAz{\TZVARptTWOx}{\TZVARptTWOy}{\TZVARptTWOz}})
--
({\TZFAx{\TZVARptTHREEx}{\TZVARptTHREEy}{\TZVARptTHREEz}},
{\TZFAy{\TZVARptTHREEx}{\TZVARptTHREEy}{\TZVARptTHREEz}}, 
{\TZFAz{\TZVARptTHREEx}{\TZVARptTHREEy}{\TZVARptTHREEz}})
--
({\TZFAx{\TZVARptFOURx}{\TZVARptFOURy}{\TZVARptFOURz}},
{\TZFAy{\TZVARptFOURx}{\TZVARptFOURy}{\TZVARptFOURz}}, 
{\TZFAz{\TZVARptFOURx}{\TZVARptFOURy}{\TZVARptFOURz}})
--cycle;
\draw[color=red, line width=0.4pt, smooth, samples=20, domain=0:\TZVARt] plot[variable=\TZVARs] ({\TZFAx{\TZVARrinter}{0}{-\TZVARrinter/\TZVARh}}, {\TZFAy{\TZVARrinter}{0}{-\TZVARrinter/\TZVARh}}, {\TZFAz{\TZVARrinter}{0}{-\TZVARrinter/\TZVARh}});
\def\TZVARtouchlen{0.605} 
\draw[color=red, dashed, line width=0.2pt] (0,0,0)--({\TZVARtouchlen*\TZVARr*cos(\TZVARt r)},{\TZVARtouchlen*\TZVARr*sin(\TZVARt r)}, {-\TZVARtouchlen*\TZVARr/\TZVARh});
}

\end{scope}
\draw (0,-1.85) node[below] {\scriptsize $\varphi=4.8$};
\end{scope}

\begin{scope}[xshift=3*\TZVARsep]

\begin{scope}[param3d]

\draw[color=black, line width=0.4pt, fill=\TZColorCone] plot[variable=\TZVARt, smooth, samples=20, domain=-pi/2:pi/2] ({\TZVARr*cos(\TZVARt r)},{\TZVARr*sin(\TZVARt r)},{-\TZVARr/\TZVARh})--(0,0,0)--cycle;
\draw[color=black, dotted, line width=0.4pt, smooth, samples=20, domain=pi/2:3*pi/2] plot[variable=\TZVARt] ({\TZVARr*cos(\TZVARt r)},{\TZVARr*sin(\TZVARt r)},{-\TZVARr/\TZVARh});

\draw[color=black, line width=0.15pt] plot[variable=\TZVARt, smooth, samples=20, domain=-pi/2:pi/2] ({\TZVARrinter*cos(\TZVARt r)},{\TZVARrinter*sin(\TZVARt r)},{-\TZVARrinter/\TZVARh});
\draw[color=black, dotted, line width=0.15pt, smooth, samples=20, domain=pi/2:3*pi/2] plot[variable=\TZVARt] ({\TZVARrinter*cos(\TZVARt r)},{\TZVARrinter*sin(\TZVARt r)},{-\TZVARrinter/\TZVARh});

\def\TZVARs{0}

\foreach \TZVARt in {(5.2)}
{
\draw[color=brown, line width=0.4pt, fill=white, opacity=0.8] 
({\TZFAx{\TZVARptONEx}{\TZVARptONEy}{\TZVARptONEz}},
{\TZFAy{\TZVARptONEx}{\TZVARptONEy}{\TZVARptONEz}}, 
{\TZFAz{\TZVARptONEx}{\TZVARptONEy}{\TZVARptONEz}})
--
({\TZFAx{\TZVARptTWOx}{\TZVARptTWOy}{\TZVARptTWOz}},
{\TZFAy{\TZVARptTWOx}{\TZVARptTWOy}{\TZVARptTWOz}}, 
{\TZFAz{\TZVARptTWOx}{\TZVARptTWOy}{\TZVARptTWOz}})
--
({\TZFAx{\TZVARptTHREEx}{\TZVARptTHREEy}{\TZVARptTHREEz}},
{\TZFAy{\TZVARptTHREEx}{\TZVARptTHREEy}{\TZVARptTHREEz}}, 
{\TZFAz{\TZVARptTHREEx}{\TZVARptTHREEy}{\TZVARptTHREEz}})
--
({\TZFAx{\TZVARptFOURx}{\TZVARptFOURy}{\TZVARptFOURz}},
{\TZFAy{\TZVARptFOURx}{\TZVARptFOURy}{\TZVARptFOURz}}, 
{\TZFAz{\TZVARptFOURx}{\TZVARptFOURy}{\TZVARptFOURz}})
--cycle;
\draw[color=red, line width=0.4pt, smooth, samples=20, domain=0:\TZVARt] plot[variable=\TZVARs] ({\TZFAx{\TZVARrinter}{0}{-\TZVARrinter/\TZVARh}}, {\TZFAy{\TZVARrinter}{0}{-\TZVARrinter/\TZVARh}}, {\TZFAz{\TZVARrinter}{0}{-\TZVARrinter/\TZVARh}});
\def\TZVARtouchlen{0.595} 
\draw[color=red, dashed, line width=0.2pt] (0,0,0)--({\TZVARtouchlen*\TZVARr*cos(\TZVARt r)},{\TZVARtouchlen*\TZVARr*sin(\TZVARt r)}, {-\TZVARtouchlen*\TZVARr/\TZVARh});
}

\end{scope}
\draw (0,-1.85) node[below] {\scriptsize $\varphi=5.2$};
\end{scope}

\begin{scope}[xshift=4*\TZVARsep]

\begin{scope}[param3d]

\draw[color=black, line width=0.4pt, fill=\TZColorCone] plot[variable=\TZVARt, smooth, samples=20, domain=-pi/2:pi/2] ({\TZVARr*cos(\TZVARt r)},{\TZVARr*sin(\TZVARt r)},{-\TZVARr/\TZVARh})--(0,0,0)--cycle;
\draw[color=black, dotted, line width=0.4pt, smooth, samples=20, domain=pi/2:3*pi/2] plot[variable=\TZVARt] ({\TZVARr*cos(\TZVARt r)},{\TZVARr*sin(\TZVARt r)},{-\TZVARr/\TZVARh});

\draw[color=black, line width=0.15pt] plot[variable=\TZVARt, smooth, samples=20, domain=-pi/2:pi/2] ({\TZVARrinter*cos(\TZVARt r)},{\TZVARrinter*sin(\TZVARt r)},{-\TZVARrinter/\TZVARh});
\draw[color=black, dotted, line width=0.15pt, smooth, samples=20, domain=pi/2:3*pi/2] plot[variable=\TZVARt] ({\TZVARrinter*cos(\TZVARt r)},{\TZVARrinter*sin(\TZVARt r)},{-\TZVARrinter/\TZVARh});

\def\TZVARs{0}

\foreach \TZVARt in {(2*pi)}
{
\draw[color=brown, line width=0.4pt, fill=white, opacity=0.8] 
({\TZFAx{\TZVARptONEx}{\TZVARptONEy}{\TZVARptONEz}},
{\TZFAy{\TZVARptONEx}{\TZVARptONEy}{\TZVARptONEz}}, 
{\TZFAz{\TZVARptONEx}{\TZVARptONEy}{\TZVARptONEz}})
--
({\TZFAx{\TZVARptTWOx}{\TZVARptTWOy}{\TZVARptTWOz}},
{\TZFAy{\TZVARptTWOx}{\TZVARptTWOy}{\TZVARptTWOz}}, 
{\TZFAz{\TZVARptTWOx}{\TZVARptTWOy}{\TZVARptTWOz}})
--
({\TZFAx{\TZVARptTHREEx}{\TZVARptTHREEy}{\TZVARptTHREEz}},
{\TZFAy{\TZVARptTHREEx}{\TZVARptTHREEy}{\TZVARptTHREEz}}, 
{\TZFAz{\TZVARptTHREEx}{\TZVARptTHREEy}{\TZVARptTHREEz}})
--
({\TZFAx{\TZVARptFOURx}{\TZVARptFOURy}{\TZVARptFOURz}},
{\TZFAy{\TZVARptFOURx}{\TZVARptFOURy}{\TZVARptFOURz}}, 
{\TZFAz{\TZVARptFOURx}{\TZVARptFOURy}{\TZVARptFOURz}})
--cycle;
\draw[color=red, line width=0.4pt, smooth, samples=20, domain=0:\TZVARt] plot[variable=\TZVARs] ({\TZFAx{\TZVARrinter}{0}{-\TZVARrinter/\TZVARh}}, {\TZFAy{\TZVARrinter}{0}{-\TZVARrinter/\TZVARh}}, {\TZFAz{\TZVARrinter}{0}{-\TZVARrinter/\TZVARh}});
\def\TZVARtouchlen{0.615} 
\draw[color=red, dashed, line width=0.2pt] (0,0,0)--({\TZVARtouchlen*\TZVARr*cos(\TZVARt r)},{\TZVARtouchlen*\TZVARr*sin(\TZVARt r)}, {-\TZVARtouchlen*\TZVARr/\TZVARh});
}

\end{scope}
\draw (0,-1.85) node[below] {\scriptsize $\varphi=2\pi$};
\end{scope}

\end{scope} 

\end{tikzpicture}
\end{figure}
}%

A rather different way to imagine the movement of the rolling tangent space uses the idea of polishing a surface. 
Imagine we have a surface $M$ in $\RR^3$, say the surface of a piece of wood, that we want to polish. 
For that we take a piece of sandpaper (a plane\slash tangent space), we apply it to the wood and we make movements to polish the wood. 
For these movements there are several different ways: we can move the sandpaper parallel to itself (translations in $\RR^3$), which is the standard way, but we can also rotate the sandpaper around its point of contact, for instance when we want to polish a small disc around this point. 
But as the piece of wood will not be flat (in general), we also have to change the plane of the sandpaper in order to follow the curvature of the piece of wood. 
In general, when polishing, one will make a combination of all three types of moves (parallel translating the sand paper, rotating it, both without changing the plane containing the sand paper, and changing this plane). 

But what should we do when we don't want to leave any scratches on surface of the wood, but simply want to \myquote{move} the piece of sandpaper along the wood \stress{without} polishing? 
Obviously we should not make parallel translations, nor rotations. 
When not making that kind of movement, (almost) all points of the sandpaper will move and their trajectories will be orthogonal to the plane containing the piece of sandpaper, as else the movement would have a component of parallel translation or rotation, which would leave a trace of polishing. 
It will be this way to visualise the rolling tangent space that will be the basis of our mathematical formulation of this movement: curves that belong at all times to the (affine) tangent space, but whose velocity vector is orthogonal to these (affine) tangent spaces. 

\medskip

In order to explain the link with parallel transport, let us start with the observation that parallel transport is basically the foundation of affine geometry! 
Let $H\subset \RR^\nu$ be an affine subspace with $E\subset\RR^\nu$ its underlying vector space.\footnote{We put $H$ inside a vector space to simplify the exposition and to avoid having to define the abstract notion of an affine space.} 
This means that to each couple of points $x,y\in H$ we can associate a vector $v=y-x\in E$, a vector that we interpret as the displacement vector of the point $x$ to the point $y$. 
We draw it as an arrow in $H$ going from $x$ to $y$. 
When doing so, we will say that $v$ is a vector attached to the (initial) point $x$. 
But this same vector $v$, representing a displacement, can be seen as the arrow going from $a\in H$ to $a+v\in H$.
\WelNietTikZpicture
{%
\begin{figure}[!ht]
\begin{tikzpicture}[scale=1, line cap=round,line join=round,x=\TZschaal, y=\TZschaal]

\newcommand\TZVARr{0.04}
\newcommand\TZVAReps{0.08}
\newcommand\puntX{(0.5,1.0)}
\newcommand\TZVARvect{(0.8,0.8)}
\newcommand\puntY{(1.3,1.8)}
\newcommand\puntA{(2.0,0.5)}
\newcommand\puntAPV{(2.8,1.3)}
\newcommand\TZFUNCf[1]{{1.4-0.5*(#1-1)*(#1-1)}}

\newcommand\TZMACarrows[2]{%
\draw[fill, color=#2] (#1,{\TZFUNCf{#1}}) circle(\TZVARr);
\draw[->, color=#2, line width=0.8pt] (#1,{\TZFUNCf{#1}})--++\TZVARvect;
}%

\begin{footnotesize}
\begin{scope}

\draw[->] (0.,-0.2) -- (0.,2.2);
\draw[->] (-0.5, 0) -- (3.5, 0);

\draw[fill] \puntX circle(\TZVARr)node[left]{$x$};
\draw[fill] \puntY circle(\TZVARr)node[right]{$y$};
\draw[fill] \puntA circle(\TZVARr)node[left]{$a$};
\draw[fill, color=gray!75] \puntAPV circle(\TZVARr);

\draw[->] \puntX++(\TZVAReps,\TZVAReps)--++(0.8-2*\TZVAReps,0.8-2*\TZVAReps) node[midway, above left]{\llap{$v$}};

\draw[->, color=gray!75] \puntA++(\TZVAReps,\TZVAReps)--++(0.8-2*\TZVAReps,0.8-2*\TZVAReps) node[midway, above left]{\llap{$v$}};

\end{scope}

\begin{scope}[xshift=7\TZschaal]

\draw[->] (0.,-0.2) -- (0.,2.2);
\draw[->] (-0.5, 0) -- (3.5, 0);

\TZMACarrows{1.45}{gray!50}
\TZMACarrows{2.05}{gray!50}

\draw[line width=0.8pt, smooth,samples=20,domain=0.2:2.8] plot[variable=\t]
({\t},{\TZFUNCf{\t}});

\draw[fill] (0.5,{\TZFUNCf{0.5}}) circle(\TZVARr)node[below]{\footnotesize $\gamma(s)$};
\draw[->, line width=0.8pt] (0.5,{\TZFUNCf{0.5}})--++\TZVARvect node[midway, above left]{\llap{$v$}};

\draw[fill] (2.5,{\TZFUNCf{2.5}}) circle(\TZVARr)node[left]{\footnotesize $\gamma(t)$};
\draw[->, line width=0.8pt] (2.5,{\TZFUNCf{2.5}})--++\TZVARvect node[midway,  right]{\ $v$ translated};

\end{scope}
\end{footnotesize}
\end{tikzpicture}
\end{figure}
}%
In doing so, we have simply \myquote{translated} the vector $v$ attached to the initial point $x$ towards the (new initial) point $a$; one says that $v$ is parallel transported from $x$ to $a$. 
Now when we add a curve $\gamma:I\to H$ and we take $x=\gamma(\tcour)$ and $a=\gamma(\tfixe)$, then one says that $v$ is parallel transported from $x=\gamma(\tcour)$ to $a=\gamma(\tfixe)$ along the curve $\gamma$ (although adding \myquote{along $\gamma$} seems to be completely fortuitous). 

But what can we do when the curve is not in an affine space but on an arbitrary submanifold $M\subset \RR^\nu$ and when we have a tangent vector $v\in T_{\gamma(\tcour)}M$, \ie, attached to the point $\gamma(\tcour)$, and we want to parallel transport $v$ to the point $\gamma(\tfixe)$? 
Because, even though we can see $v$ as a vector in the underlying (affine) vector space $\RR^\nu$, it is highly unlikely that the parallel transported vector to the point $\gamma(t)$ will be tangent to $M$ at $\gamma(\tfixe)$, as in general the tangent spaces vary from point to point. 
The rolling tangent space gives us a possible answer to this question. 

Let us go back to our idea of a person inside a ball $M\subset \RR^3$. 
In order to distinguish between points\slash vectors related to the ball $M$ and points\slash vectors related to the (flat) ground (the affine tangent space) represented by $\RR^2$, we will prefix the latter by the letter $\ground$ for \myquote{ground.} 
When at time $\tcour$ his feet are at the point $\gamma(\tcour)\in M$, this point (necessarily) touches the ground in a point (the same point) denoted by $\ground{\gamma(\tcour)}\in \RR^2$. 
A tangent vector $v\in T_{\gamma(\tcour)}M$ thus is at the same time a vector $\ground v$ on the ground, a vector attached to the point $\ground{\gamma(\tcour)}$. 
When the person inside the ball moves along the curve $\gamma$ painted on the ball towards the point $\gamma(\tfixe)\in M$ at time $\tfixe$, the ball rolls (without sliding) on the ground and thus the point $\gamma(\tcour)$ and the tangent vector $v\in T_{\gamma(\tcour)}M$ move with respect to a fixed spectator of this game, whereas the point $\ground\gamma(\tcour)$ and the vector $\ground v$ attached to this point remain fixed on the ground. 
In this way we obtain the trace curve on the ground, but we also have the vector $\ground v$ attached to $\ground\gamma(\tcour)$, the imprint of the tangent vector $v\in T_{\gamma(\tcour)}M$. 
\WelNietTikZpicture
{
\begin{figure}[!ht]
\begin{tikzpicture}[scale=1, line cap=round,line join=round,x=\TZschaal, y=\TZschaal]
\def\TZemme{0.2} 
\def\TZhoek{30} 
\def\TZVARr{0.03}

\begin{scope}[param3d, xshift=1.5\TZschaal]
\def\TZVARrectsize{2.2}
\draw[color=green!70, line width=0.4pt] ({-\TZVARrectsize}, {-\TZVARrectsize}, -1)--({\TZVARrectsize}, {-\TZVARrectsize}, -1)--({\TZVARrectsize}, {\TZVARrectsize}, -1)--({-\TZVARrectsize}, {\TZVARrectsize}, -1)--cycle;
\end{scope}

\draw[color=brown, line width=0.5pt, fill=white, xshift=3\TZschaal, yshift=0.05\TZschaal] (0,0) circle (1);

\draw[fill, xshift=3\TZschaal, yshift=0.05\TZschaal] (0,-1)circle(\TZVARr);
\draw[xshift=3\TZschaal, yshift=0.05\TZschaal] (0.2,-1) node[above]{\tiny $\gamma(\tcour)$\kern1em};

\draw[fill, xshift=3\TZschaal, yshift=0.05\TZschaal] (-0.35,0.76)circle(\TZVARr);
\draw[xshift=3\TZschaal, yshift=0.05\TZschaal] (-0.35,0.76) node[right]{\tiny $\gamma(\tfixe)$};

\begin{scope}[param3d, xshift=3\TZschaal, yshift=0.05\TZschaal]

\def\TZVARphione{(-60)}\def\TZVARphitwo{(120)}
\draw[color=brown, line width=0.5pt, smooth, samples=20, domain=\TZVARphione:\TZVARphitwo] plot[variable=\t] ({cos(\t)}, {sin(\t)},{0});
\draw[color=brown, dotted, line width=0.5pt, smooth, samples=20, domain=\TZVARphitwo:\TZVARphione+360] plot[variable=\t] ({cos(\t)}, {sin(\t)},{0});

\def\TZVARphi{-20}
\draw[color=red, line width=0.5pt, smooth, samples=20, domain=-90:-20] plot[variable=\t] 
({cos(\TZVARphi)*cos(\t)},{sin(\TZVARphi)*cos(\t)},{sin(\t)}) to[out=105, in=-130](-1,-1.55,0.2)to[out=50, in=180](-1,-1,0.5);

\draw[->] (0,0,-1)--(0.7,0.7,-1) node[right]{\tiny $v$};
\end{scope}
\begin{scope}[param3d]

\draw[color=blue, dotted, line width=0.5pt, smooth, samples=20, domain=0:2.48] plot[variable=\t] ({1-cos(2*\t r)-\t},{\t},{-1});

\end{scope}

\begin{scope}[xshift=7.5\TZschaal]

\begin{scope}[param3d, xshift=1.5\TZschaal]
\def\TZVARrectsize{2.2}
\draw[color=green!70, line width=0.4pt] ({-\TZVARrectsize}, {-\TZVARrectsize}, -1)--({\TZVARrectsize}, {-\TZVARrectsize}, -1)--({\TZVARrectsize}, {\TZVARrectsize}, -1)--({-\TZVARrectsize}, {\TZVARrectsize}, -1)--cycle;
\end{scope}

\draw[color=brown, line width=0.5pt, fill=white] (0,0) circle (1);

\begin{scope}[param3d]

\def\TZVARphione{(-60)}\def\TZVARphitwo{(120)}
\draw[color=brown, line width=0.5pt, smooth, samples=20, domain=\TZVARphione:\TZVARphitwo] plot[variable=\t] ({cos(\t)}, {sin(\t)},{0});
\draw[color=brown, dotted, line width=0.5pt, smooth, samples=20, domain=\TZVARphitwo:\TZVARphione+360] plot[variable=\t] ({cos(\t)}, {sin(\t)},{0});

\draw[color=blue, line width=0.5pt, smooth, samples=20, domain=0:2.48] plot[variable=\t] ({1-cos(2*\t r)-\t},{\t},{-1});

\draw[color=red, line width=0.5pt, smooth, samples=20, domain=0:0.9] plot[variable=\t] ({1-cos(2*\t r)-\t},{\t},{sqrt(2*(1-\t)*(cos(2*\t r)+\t)-(cos(2*\t r)*cos(2*\t r)))*(-1)}) to[in=-90, out=86](1,1,1);

\draw[->, color=gray] (0,0,-1)--(0.7,0.7,-1);

\draw[->, xshift=3\TZschaal, yshift=0.05\TZschaal] (0,0,-1)--(0.7,0.7,-1) node[right]{\tiny $\ground v$};

\draw[->] (0,0.425,0.77)--++(0,-0.4,-0.1) node[left]{\tiny $v$};

\end{scope}

\draw[fill] (0,-1)circle(\TZVARr);
\draw (0.1,-0.9) node{\tiny $\gamma(\tfixe)$\vrule width0pt depth1.5em\kern1em};
\draw[color=black!90, dotted] (0.35,-1.2) to[out=-45, in=180] (0.35,-1.7) node[right]{\tiny $\ground v$ parallel transported};

\draw[fill] (0.366,0.71)circle(\TZVARr);
\draw (0.4,0.61) node[above right]{\tiny $\gamma(\tcour)$};

\draw[fill, xshift=3\TZschaal, yshift=0.05\TZschaal] (0,-1)circle(\TZVARr);
\draw[xshift=3\TZschaal, yshift=0.05\TZschaal] (0.2,-1) node[above]{\tiny $\ground \gamma(\tcour)$};

\end{scope}
\end{tikzpicture}

\end{figure}
}%
We then can \myquote{parallel transport} this trace vector $\ground v$ on the ground from the initial point $\ground{\gamma(\tcour)}$ to the new contact point on the ground $\gamma(\tfixe) \equiv \ground\gamma(\tfixe)$. 
This parallel transported vector thus will be simultaneously a vector on the ground, \ie, the new\slash old affine tangent space, as well as a tangent vector to $M$ at the point $\gamma(\tfixe)$. 
In this way we have obtained a tangent vector at $T_{\gamma(\tfixe)}M$ associated to the initial tangent vector $v\in T_{\gamma(\tcour)}M$ by means of parallel transport in the rolling tangent space. 
And as we will show, it is exactly this procedure that is called parallel transport of a tangent vector along a curve.

\section{An overview of the mathematical results}
\label{SECTIONSURPLANTANGENTBASCULANTEVOLUTE}

Starting here, the six symbols $I$, $\nu$, $M$, $n$, $\gamma$ and $\pi$ will have a fixed meaning: $I$ will be an interval in $\RR$, $\nu\in\NN^*$ the dimension of the surrounding vector space $\RR^\nu$, $M\subset \RR^\nu$ will be a submanifold\footnote{It is on purpose that we do not give the definition of a submanifold. On the one hand, this definition is not really important for us and on the other hand, giving it would require far too many explanations more or less without a purpose. The only really important point for us is that the tangent spaces $T_{\gamma(\tcour)}M$ depend differentiably on $\tcour$.} of dimension $n$ and class $C^2$, $\gamma:I\to M$ will be a curve of class $C^2$ on $M$ and $\pi:\RR^\nu\to T_{\gamma(\tcour)}M$ will denote the orthogonal projection of $\RR^\nu$ on the tangent space  to $M$ at the point $\gamma(\tcour)$. 
Additionally, in order to avoid angles and other undesirable situations, we will suppose that $\gamma$ is a regular curve, \ie, $\gamma'(\tcour)\neq0$ for all $\tcour\in I$. 
We will also use $\inprod vw$ for the standard scalar product on $\RR^\nu$ and $\norme v = \sqrt{\inprod vv}$ for the associated norm. 

We define the family of affine subspaces $H_\tcour \subset \RR^\nu$, $\tcour\in I$  associated to $\gamma$ by
\begin{moneq}
H_\tcour = \gamma(\tcour)+T_{\gamma(\tcour)}M
\mapob.
\end{moneq}
The vector space underlying the affine space $H_\tcour$ is the tangent space $T_{\gamma(\tcour)}M$. 
As explained in \recals{SECTIONSURPLANTANGENTBASCULANT}, in order to find the movement of the rolling tangent space\footnote{The reader who thinks that the drawing below of this movement of the tangent space is familiar is quite right: the movement of the rolling tangent space is a direct generalisation of the involute of a curve.} we will be looking for curves $x:I\to \RR^\nu$ that satisfy the following two conditions:
\begin{enumerate}[\ (\labelRTS1)\,]
\item\label{RTS1}
$\forall t\in I: x(\tcour)\in H_\tcour$ and 

\item\label{RTS2}
$\forall \tcour\in I: \inprod{x'(\tcour)}{T_{\gamma(\tcour)}M}=0$, \ie, $x'(\tcour)$ is orthogonal to $H_\tcour$.

\end{enumerate}
\WelNietTikZpicture
{
\begin{figure}[!ht]
\begin{tikzpicture}[scale=1, line cap=round,line join=round,x=\TZschaal, y=\TZschaal]

\def\TZFUNCx#1{cos((#1) r)+((#1)-\TZVARc*pi/5+pi/2)*sin((#1) r)}
\def\TZFUNCy#1{sin((#1) r) + (\TZVARc*pi/5-pi/2-(#1))*cos((#1) r)}

\draw (0,1) node[below]{\scriptsize$\gamma(\tcour)$};
\draw ({cos((2*pi/5-pi/2) r)},{sin((2*pi/5-pi/2) r)}) node[left]{\scriptsize$\gamma(t')$};
\draw ({cos((-pi/5-pi/2) r)},{sin((-pi/5-pi/2) r)}) node[above ]{\scriptsize$\kern1em\gamma(\tfixe)$};
\def\TZVARc{5}
\draw[color=black] ({\TZFUNCx{-pi/2-pi/5}},{\TZFUNCy{-pi/2-pi/5}}) node[right]{\scriptsize$(\TraceAff\gamma)(\tcour)\vrule width0pt depth1.5ex$};
\def\TZVARc{2}
\draw[color=black] ({\TZFUNCx{-pi/2-pi/5}},{\TZFUNCy{-pi/2-pi/5}}) node[below]{\scriptsize\kern0.4em$(\TraceAff\gamma)(t')$};
\def\TZVARc{-4}
\draw ({0.2+\TZFUNCx{-pi/2-pi/5}},{\TZFUNCy{-pi/2-pi/5}}) node[left]{\scriptsize$H_\tfixe$};
\draw ({3.9},{1}) node[right]{\scriptsize$H_\tcour$};
\def\TZVARc{5}
\draw ({0.1},{\TZFUNCy{-pi/2-pi/5}}) node[left]{\scriptsize$H_{t'}$};
\draw[color=red, dotted] (3.5,0)--(3.99,0) node[right]{\scriptsize : the movement of the};
\draw[color=red] (3.35,-0.3) node[right]{\scriptsize \kern3.5em rolling tangent space};

\clip (-2,-3.2) rectangle (4,1.6);

\draw[color=black, line width=1.5pt] (0,0) circle (1);

\foreach \TZVARrot in {-1,...,5}
{
\draw[color=black, line width=0.2pt, rotate=deg(\TZVARrot*pi/5), shift={({-\TZVARrot*pi/5-1},{-1})}] (-10,0)--(10,0);
\draw[fill, color=black, line width=0.2pt, rotate=deg(\TZVARrot*pi/5), shift={({0},{-1})}] (0,0) circle(0.05);
}

\foreach\TZVARc in {2,5}
{
\draw[color=red,dotted, line width=0.4pt, smooth,samples=20,domain=-pi/2-pi/5:{\TZVARc*pi/5-pi/2}] plot[variable=\t] ({\TZFUNCx{\t}},{\TZFUNCy{\t}});
\draw[color=red,dotted, line width=0.4pt, smooth,samples=20,domain={\TZVARc*pi/5-pi/2}:{-pi/2+pi}] plot[variable=\t] ({\TZFUNCx{\t}},{\TZFUNCy{\t}});
\foreach \x in {-1,...,5}
{
\draw[fill, color=red] ({\TZFUNCx{\x*pi/5-pi/2}},{\TZFUNCy{\x*pi/5-pi/2}}) circle({0.035*(\TZVARc>\x)});
\draw[fill, color=red] ({\TZFUNCx{\x*pi/5-pi/2}},{\TZFUNCy{\x*pi/5-pi/2}}) circle({0.035*(\TZVARc<\x)});
}
}

\end{tikzpicture}
\end{figure}%
}%
In \recalt{ExpressionExpliciteBasulementTangentEtTransportParallel} we will show that for all $\tfixe\in I$ and all $x_o\in H_\tfixe$ there exists a unique curve $x_{\tfixe,x_o}:I\to \RR^\nu$ that satisfies these two conditions as well as the initial condition $x_{\tfixe,x_o}(\tfixe) = x_o$. 
This allows us to define the \stresd{movement of the rolling tangent space $\esptanbas$} as the family of maps  $\esptanbas(\tcour,\tfixe):H_\tfixe\to H_\tcour$ (for $\tfixe,\tcour\in I$) by
\begin{moneq}
\esptanbas(\tcour,\tfixe)(x_o) = x_{\tfixe,x_o}(\tcour)
\mapob.
\end{moneq}
Once we have defined the movement of the rolling tangent space, we can define the \stresd{trace curve in the affine space $H_\tfixe$} as the curve $\TraceAff\gamma\NoTraceAff{\gammah_\tfixe}:I\to H_\tfixe$ that has the property that $(\TraceAff\gamma\NoTraceAff{\gammah_\tfixe})(\tcour)$ is the point in $H_\tfixe$ obtained by the movement of the rolling tangent space starting at the point $\gamma(\tcour)\in H_\tcour$. 
More precisely:
\begin{moneq}
(\TraceAff\gamma\NoTraceAff{\gammah_\tfixe})(\tcour) = \esptanbas(\tfixe,\tcour)\bigl(\gamma(\tcour)\bigr)
\mapob.
\end{moneq}

In order to prove the claims made in the introduction, we will need some intermediate results. 
The first of these intermediate results is that uniqueness of the curves $x_{\tfixe,x_o}$ implies that we have
\begin{moneq}[LoiDeGroupeRTB]
\forall s,t,u\in I
\quad:\quad
\esptanbas(s,t)\scirc \esptanbas(t,u)= \esptanbas(s,u)
\mapob.
\end{moneq}
This in turn implies in particular that $\esptanbas(t,s)$ is bijective with $\esptanbas(t,s)\mo = \esptanbas(s,t)$. 
A second intermediate result is that for fixed $s,t\in I$, the map $\esptanbas(\tfixe,\tcour):H_\tcour\to H_\tfixe$ is an isometry \recalt{ETBetTPEstIsometrie}, hence an affine map \recalt{AppIsometrieDoncAffine}. 
The linear map underlying the affine map $\esptanbas(t,s)$ is called the \stresd{map of parallel transport} (a detailed explanation for this name will follow) and is denoted as $\transportparallel(\tfixe,\tcour):T_{\gamma(\tcour)}M\to T_{\gamma(\tfixe)}M$. 
It is given by the formula
\begin{moneq}[TransportParallelParDifferenceEspTanBas]
\begin{aligned}
\transportparallel(\tfixe,\tcour)(x-y)
&
=\esptanbas(\tfixe,\tcour)(x) -\esptanbas(\tfixe,\tcour)(y)
\\
\text{or}\qquad
\transportparallel(\tfixe,\tcour)(v)
&
=\esptanbas(\tfixe,\tcour)(x_o+v) -\esptanbas(\tfixe,\tcour)(x_o)
\mapob,
\end{aligned}
\end{moneq}
for an aribtrary $x_o\in H_\tcour$. 
And because $\esptanbas(\tfixe,\tcour)$ is an isometry, $\transportparallel(\tfixe,\tcour)$ is too. 

We now have all the necessary ingredients to translate our intuitive idea into formul\ae. 
As the tangent space $T_{\gamma(\tcour)}M$ is the vector space underlying the affine space $H_\tcour$, a tangent vector $v\in T_{\gamma(\tcour)}M$ can be realised as a couple of points $x,y\in H_\tcour$ satisfying $v=y-x$, in which case $v$ would be a vector attached to the (initial) point $x$. 
The movement of the rolling tangent space $\esptanbas(\tfixe,\tcour):H_\tcour\to H_t$ sends (transports) the two points $x,y\in H_\tcour$ to two new points $x',y'\in H_\tfixe$, and these two points can be interpreted as a vector $y'-x'=v' \in T_{\gamma(\tfixe)}M$ (the vector space underlying the affine space $H_\tfixe$) attached to the (initial) point $x' = \esptanbas(\tfixe,\tcour)(x)$. 
Knowing that $\esptanbas(\tfixe,\tcour)$ is an affine map whose underlying linear map is $\transportparallel(\tfixe,\tcour)$, this new vector $v'$ attached to $x'$ is given by $v'= \transportparallel(\tfixe,\tcour)(v)$. 
Taking $x=\gamma(\tcour)$, which is quite natural as $v$ belongs to the tangent space to $M$ at the point $\gamma(\tcour)$, the image vector $v'=\transportparallel(\tfixe,\tcour)(v)\in T_{\gamma(\tfixe)}M$ will be attached to the image point $x'=\esptanbas(\tfixe,\tcour)\bigl(\gamma(\tcour)\bigr) = (\TraceAff\gamma\NoTraceAff{\gammah_\tfixe})(\tcour)$, a point of the trace curve $\TraceAff\gamma\NoTraceAff{\gammah_\tfixe}$ in $H_\tfixe$. 
But $H_\tfixe$ is an affine space, so we can \myquote{attach} the vector $v'$ also to the point $\gamma(\tfixe)\in H_\tfixe$. 
In this way we have \myquote{transported} the vector $v'$ from (being attached to) $(\TraceAff\gamma\NoTraceAff{\gammah_\tfixe})(\tcour)$ to (being attached to) $\gamma(\tfixe)$. 

To summarise: starting with a tangent vector $v\in T_{\gamma(\tcour)}M$ at $\gamma(\tcour)$, we apply the movement of the rolling tangent space to obtain a vector in $H_\tfixe$ attached to the point $(\TraceAff\gamma\NoTraceAff{\gammah_\tfixe})(\tcour)$, which we then \myquote{parallel transport} to attach it to the point $\gamma(\tfixe)\in H_\tfixe$. 
It is this procedure that is called the \stresd{parallel transport of $v\in T_{\gamma(\tcour)}M$ to $T_{\gamma(\tfixe)}M$ along the curve $\gamma$}. 
In this general situation it is important to add \myquote{along $\gamma$,} as in general the results of this parallel transport of $v\in T_{\gamma(\tcour)}M$ to $v'\in T_{\gamma(\tfixe)}M$ will depend upon the curve. 
We thus have {shown} that our intuitive idea of parallel transport can be formulated by formul\ae{} as
\begin{gather}
v'\in T_{\gamma(\tfixe)}M \text{ is the parallel transport of } v\in T_{\gamma(\tcour)}M \text{ to } T_{\gamma(\tfixe)}M
\notag
\\[-1\jot]
\text{or}
\notag
\\[-1\jot]
v\in T_{\gamma(\tcour)}M\text{ is parallel transported to }v'\in T_{\gamma(\tfixe)}M 
\label{NotreDefinitionTransportParallele}
\\
\Updownarrow
\notag
\\
v' = \transportparallel(\tfixe,\tcour)(v)
\mapob,
\notag
\end{gather}
which explains the name \stresd{parallel transport} for the map $\transportparallel(\tfixe,\tcour)$. 
But beware: the two tangent spaces $T_{\gamma(\tcour)}M$ and $T_{\gamma(\tfixe)}M$ are different and there is no reason to think that the vectors $v$ and $v'$ are parallel, and even less that they are equal. 
We can only be sure that they have the same length: as $\transportparallel(\tfixe,\tcour)$ is an isometry, the vectors $v$ and $v' = \transportparallel(\tfixe,\tcour)(v)$ have the same length. 

Let us now look at the notion of a field $v$ of tangent vectors along $\gamma$, meaning that $v:I\to \RR^\nu$ is a map satisfying the condition that for all $\tcour\in I$ we have $v(\tcour)\in T_{\gamma(\tcour)}M$.
(Note that we do not require that $v(\tcour)$ is tangent to the curve $\gamma$, only that it is tangent to $M$.) 
We have argued above that $\transportparallel(s,t)\bigl(v(\tcour)\bigr)$ is (can be interpreted as) a vector attached to the point $(\TraceAff\gamma\NoTraceAff{\gammah_\tfixe})(\tcour)$ of the trace curve  inside the affine space $H_\tfixe$. 
In this way we obtain a new field of tangent vectors $\TraceAff v : I\to T_{\gamma(\tfixe)}M$ (in this case tangent to $H_\tfixe$) along the trace curve $\TraceAff\gamma\NoTraceAff{\gammah_\tfixe}$ defined by
\begin{moneq}[DefinitionChampTraceLeLongCourbe]
(\TraceAff v)(\tcour) = \transportparallel(\tfixe,\tcour)\bigl(v(\tcour)\bigr)
\mapob.
\end{moneq}
And on says that $v$ is a \stresd{field of parallel (tangent) vectors} along $\gamma$ if this trace field $\TraceAff v$ is constant (necessarily equal to $v(\tfixe)$). 
Said differently --- when we displace the point of attachment of the vector $(\TraceAff v)(\tcour)$ from the point $(\TraceAff\gamma)(\tcour)$ to the point $\gamma(\tfixe)$ --- $v$ is a field of parallel (tangent) vectors if for all $\tcour\in I$ the parallel transport \recalf{NotreDefinitionTransportParallele} of $v(\tcour)\in T_{\gamma(\tcour)}M$ to $T_{\gamma(\tfixe)}M$ yields the vector $v(\tfixe)$. 

In order to make the link with the usual definition of a field $v$ of parallel vectors, we note that the latter uses the definition of the covariant derivative of $v$. 
We start with the observation that, even if $v(\tcour)\in T_{\gamma(\tcour)}M$ for all $\tcour\in I$, there is no reason to believe that $v'(\tcour)\in T_{\gamma(\tcour)}M$. 
In order to stay inside this tangent space one thus defines the \stresd{covariant derivative of $v$}, noted as $Dv/\extder\tcour$, as the orthogonal projection of the ordinary derivative $v'(\tcour)$ on $T_{\gamma(\tcour)}M$ (\cite[10.4.7 p376]{BG88}, \cite[\S4-4 p238]{DoC76}, \cite[\S7.4 p171]{Pressley..2010}, \cite[5.37 p280]{Tapp..2016}), \ie,
\begin{moneq}[DefinitionUsuelleDeriveeCovariante]
\frac{Dv}{\extder \tcour}(\tcour) \oversettext{\ def.\ }\to= \pi\bigl(v'(\tcour)\bigr)
\mapob,
\end{moneq}
which thus produces a new field of tangent vectors $Dv/\extder \tcour$  along $\gamma$. 
One then says that $v$ is a \stresd{field of parallel vectors} if its covariant derivative is (identically) zero. 
To show that this definition is equivalent to ours we use \recalt{TransportparalleleDeGammahpEtGammahsec}, which tells us that the ordinary derivative of the trace field $\TraceAff v\NoTraceAff{\vh_\tfixe}$ and the ordinary derivative of $v$ are related by the formula 
\begin{moneq}[LienDeriveeCovarianteEtDeriveeOrdinaireTraceChamp]
\transportparallel(\tcour,\tfixe)\bigl((\TraceAff v)'(\tcour)\NoTraceAff{(\vh_\tfixe)'(\tcour)}\bigr) 
\oversettext{\ \recalt{TransportparalleleDeGammahpEtGammahsec}\ }\to= 
\pi\bigl(v'(\tcour)\bigr) = \frac{Dv}{\extder \tcour}(\tcour)
\mapob.
\end{moneq}
As the field of tangent vectors $Dv/\extder \tcour$ also has a trace field $\TraceAff(Dv/\extder \tcour)$ on $H_\tfixe$, we can rewrite \recalf{LienDeriveeCovarianteEtDeriveeOrdinaireTraceChamp} as
\begin{moneq}[TraceDeriveeCovarianteEstDeriveeOrdinaireTrace]
\boxed{\vrule width0pt height4ex depth2.6ex
\biggl(\TraceAff\hbiggl(\frac{Dv}{\extder \tcour}\hbiggr)\biggr)(\tcour) 
= 
\frac{\extder (\TraceAff v)}{\extder \tcour}(\tcour)
}%
\equiv
(\TraceAff v)'(\tcour)
\mapob.
\end{moneq}
And because $(\TraceAff v)(\tcour)$ is zero if and only if $v(\tcour)$ is zero \recalf{DefinitionChampTraceLeLongCourbe}, it follows that $v$ is a field of parallel vectors in the ordinary sense (covariant derivative zero) if and only if the trace field $\TraceAff v$ is constant (ordinary derivative zero), which is our definition. 
Except that with the standard definition one looses the idea that parallel transport corresponds to the displacement of a vector from one point of attachement to another in an affine space (the fixed affine tangent space $H_\tfixe$). 
And one also looses the fact that the covariant derivative $Dv/\extder \tcour$ corresponds (va the isometry $\transportparallel(\tfixe,\tcour)$) to the ordinary derivative $(\TraceAff v)'\NoTraceAff{\vh_\tfixe'}$ of the trace field $\TraceAff v\NoTraceAff{\vh_\tfixe}$ inside the afine space $H_\tfixe$ \recalf{DefinitionChampTraceLeLongCourbe}, \recalf{TraceDeriveeCovarianteEstDeriveeOrdinaireTrace}. 

There is yet another way to link the standard definition of a field of parallel vectors with ours: the covariant derivative of a field $v$ is zero if (and only if) $v'$ is orthogonal to the tangent space. 
Said differently, $v:I\to \RR^\nu$ is a field of parallel (tangent) vectors along $\gamma$ if (and only if) it satisfies the two conditions $v(\tcour)\in T_{\gamma(\tcour)}M$ and $\inprod{v'(\tcour)}{T_{\gamma(\tcour)}M}=0$.\footnote{These are the conditions one finds in \cite[10.4.6, p376]{BG88} as definition of a field of parallel vectors.} 
These are the conditions of the rolling tangent space in which the affine space $H_\tcour$ is replaced by the vector space $T_{\gamma(\tcour)}M$.

\section{The link with the approach by Levi-Civita}

As said in the introduction, the idea of the rolling tangent space already appears in \cite{Persico..1921}, just 4 years after the introduction of the concept of parallel transport by Levi-Civita in 1917 \cite{LeviCivita..1917}. 
And Levi-Civita incorporated the description made by Persico in his book \cite{LeviCivita..1927}. 
In order to make the link with our description, we imagine being on a surface $M$ (of dimension $2$ in our $3$-dimensional space $\RR^3$), a surface on which we have a curve $\gamma:I\to M$. 
In \cite{Persico..1921} Persico speaks about rolling the surface on a plane, which is exactly our image of a person inside a big ball who walks on the ground. 
And his description of parallel transport is also exactly what we described. 
Except that \cite{Persico..1921} is just two pages of text without any formula, whereas we can give explicit formul\ae{} for this (rolling) movement. 
For that we recall that $\esptanbas(\tfixe,\tcour):H_\tcour\to H_\tfixe$  is an isometry and that $H_\tcour$ and $H_\tfixe$ are two affine subspaces of $\RR^3$ of dimension $2$, \ie, of codimension $1$. 
It follows that there exists a unique \myquote{positive} isometry $\esptanbast(\tfixe,\tcour):\RR^3\to\RR^3$ consisting of a rotation (not a mirror symmetry)
 and a translation whose restriction to $H_\tcour$ is $\esptanbas(\tfixe,\tcour)$. 
\WelNietTikZpicture
{%
\begin{figure}[!ht]
\begin{tikzpicture}[scale=0.75, line cap=round,line join=round,x=\TZschaal, y=\TZschaal]
\def\TZVAReps{0.06}

\begin{scope}[xshift=0\TZschaal]
\begin{scope}[rotate=30]
\draw[fill, color=red](0,0) circle[radius=\TZVAReps];
\draw[fill, color=red] (-1-\TZVAReps,-\TZVAReps)rectangle(-1+\TZVAReps, \TZVAReps);

\draw[line width=1pt] (-1,0)--(1,0)node[right]{\scriptsize$H_\tcour$};

\end{scope}

\draw[->, line width=0.75pt] (1.2,-0.5) to[out=-30, in=-150] ++(1.5,0);
\draw(1.95,-0.6)node[ below]{\scriptsize$\esptanbas(\tfixe,\tcour)$};

\begin{scope}[xshift=4\TZschaal, rotate=-20]
\draw[fill, color=red](0,0) circle[radius=\TZVAReps];
\draw[fill, color=red] (-1-\TZVAReps,-\TZVAReps)rectangle(-1+\TZVAReps, \TZVAReps);

\draw[line width=1pt] (-1,0)--(1,0)node[right]{\scriptsize$H_\tfixe$};

\end{scope}
\end{scope}

\draw (7.7,-0.5)node[right]{\scriptsize$\Longrightarrow$};

\begin{scope}[xshift=12\TZschaal]
\begin{scope}[rotate=30]
\draw[fill, color=red](0,0) circle[radius=\TZVAReps];
\draw[fill, color=red] (-1-\TZVAReps,-\TZVAReps)rectangle(-1+\TZVAReps, \TZVAReps);
\draw[color=red] (-0.8,0.5) circle[radius=\TZVAReps];
\draw[line width=1pt] (-1,0)--(1,0);
\foreach \TZVARk in {-5,-4,..., 5} \draw[line width=0.4pt] (-0.8,{\TZVARk/10})--(0.8,{\TZVARk/10});
\end{scope}

\draw[->, line width=0.75pt] (1.2,-0.5) to[out=-30, in=-150] ++(1.5,0);
\draw(1.95,-0.6)node[ below]{\scriptsize$\widetilde\esptanbas(\tfixe,\tcour)$};

\begin{scope}[xshift=4\TZschaal, rotate=-20]
\draw[fill, color=red](0,0) circle[radius=\TZVAReps];
\draw[fill, color=red] (-1-\TZVAReps,-\TZVAReps)rectangle(-1+\TZVAReps, \TZVAReps);
\draw[color=red] (-0.8,0.5) circle[radius=\TZVAReps];
\draw[line width=1pt] (-1,0)--(1,0);
\foreach \TZVARk in {-5,-4,..., 5} \draw[line width=0.4pt] (-0.8,{\TZVARk/10})--(0.8,{\TZVARk/10});
\end{scope}
\end{scope}

\end{tikzpicture}
\end{figure}
}%
With these isometries we can define, for fixed $\tfixe$, the family of surfaces $M_\tcour=\esptanbast(\tfixe,\tcour)(M)\subset\RR^3$, surfaces that are thus all isometrically equivalent to $M$. 
And moreover, each $M_\tcour$ is tangent to $H_\tfixe$ at $\esptanbas(\tfixe,\tcour)\bigl(\gamma(\tcour)\bigr) = (\TraceAff\gamma)(\tcour)\NoTraceAff{\gammah_\tfixe(\tcour)}\in H_\tfixe$. 
This family thus represents the movement of rolling $M$ on the (affine) plane $H_\tfixe$ along the trace curve $\TraceAff\gamma\NoTraceAff{\gammah_t}$ as promised. 

Instead of giving formul\ae{} for this movement, Persico passes through the intermediate of a developable ruled surface, the envelope of the tangent planes $T_{\gamma(\tcour)}M$ to the curve $\gamma$. 
In order to find a description of this surface (in a very naive way), we imitate the naive idea for the envelope of a family of lines in $\RR^2$ in which we take the set of intersection points of nearby lines and we take the limit as the lines a closer and closer. 
\WelNietTikZpicture
{%
\begin{figure}[!ht]
\begin{tikzpicture}[xscale=2, line cap=round,line join=round,x=\TZschaal, y=\TZschaal]
\def\TZVARlong{1}
\def\TZVAReps{0.03}
\def\TZVARinter{0.22}

\clip (-2,0.2) rectangle (2,1.3);
\foreach\th in {-75,-50,...,75}
{
\draw ({sin(\th)+\TZVARlong*cos(\th)},{cos(\th)-\TZVARlong*sin(\th)})--++({-2*\TZVARlong*cos(\th)},{2*\TZVARlong*sin(\th)});
\draw[fill, color=black] ({sin(\th)+\TZVARinter*cos(\th)},{cos(\th)-\TZVARinter*sin(\th)}) circle[radius=\TZVAReps];
}

\end{tikzpicture}
\end{figure}
}%
For the envelope of the affine spaces $H_\tcour=\gamma(\tcour)+T_{\gamma(\tcour)}M$ we do the same: we take the set of intersection lines of nearby planes $H_\tcour$ and $H_{\tcour+\delta \tcour}$ and we take the limit $\delta \tcour\to 0$. 
In order to describe the plane $H_\tcour$, it suffices to give a unit vector $\normal(\tcour)$ orthogonal to $T_{\gamma(\tcour)}M$. 
(There exist only two different differentiable choices (in $\tcour$) for $\normal(\tcour)$, which differ by a global sign.) 
With these unit normal vectors, the direction of the intersection line then is $\normal(\tcour)\wedge \normal(\tcour+\delta \tcour)$. 
\WelNietTikZpicture
{%
\begin{figure}[!ht]
\begin{tikzpicture}[scale=0.9, line cap=round,line join=round,x=\TZschaal, y=\TZschaal]
\def\TZVAReps{0.07}
\def\TZVARsize{3}
\def\TZVARhoog{0.48}
\def\TZVARaxe{0.4}
\def\TZVARfracr{0.37}
\def\TZVARfracl{1}
\def\TZVARfracn{0.11}
\def\TZVARax{1.3}\def\TZVARay{-0.21}

\clip (1.9-2*\TZVARsize,-0.8) rectangle (5+2*\TZVARsize,1.6);

\draw (-\TZVARsize,\TZVARhoog)node[left]{\scriptsize$H_{\tcour+\delta \tcour}$}--++(2*\TZVARsize,-2*\TZVARhoog)--++(1,0.75)--++(-\TZVARfracr*\TZVARsize,\TZVARhoog*\TZVARfracr);
\draw (-\TZVARsize,-\TZVARhoog)--++(2*\TZVARsize,2*\TZVARhoog)--++(1,0.75)node[right]{\scriptsize$H_{\tcour}$}--++(-\TZVARsize,-\TZVARhoog)--++(-\TZVARsize,\TZVARhoog)--++(-1,-0.75);
\draw (0,0)--(1,0.75);
\draw[->] (1,0.75)--++(\TZVARaxe,0.75*\TZVARaxe)node[above]{\scriptsize$\normal(\tcour)\wedge\normal(\tcour+\delta \tcour)$};
\draw(-\TZVARsize,-\TZVARhoog)--++(\TZVARfracl,0.75*\TZVARfracl);

\draw (0,-8.11)circle[radius=8];
\draw[fill] (-\TZVARax,\TZVARay)circle[radius=\TZVAReps/2]node[below]{\scriptsize$\gamma(\tcour)$};
\draw[fill] (\TZVARax,\TZVARay)circle[radius=\TZVAReps/2]node[below]{\scriptsize$\gamma(\tcour+\delta \tcour)$};
\draw[->] (-\TZVARax,\TZVARay)--++(-\TZVARhoog*\TZVARfracn,\TZVARsize*\TZVARfracn);
\draw[->] (\TZVARax,\TZVARay)--++(\TZVARhoog*\TZVARfracn,\TZVARsize*\TZVARfracn);
\draw (-\TZVARax,\TZVARay+0.175)node[left]{\scriptsize$\normal(\tcour)$};
\draw (\TZVARax,\TZVARay+0.175)node[right]{\scriptsize$\normal(\tcour+\delta \tcour)$};

\begin{scope}[xshift=8\TZschaal]
\def\TZVARxc{0}\def\TZVARyc{-1.4}\def\TZVARray{1.5}
\def\TZFONext#1#2{%
({\TZVARray*cos(#1)/cos(#2)},{\TZVARyc+\TZVARray*sin(#1)/cos(#2)})
}
\def\TZFONborder#1#2{%
({#2*cos(#1)},{#2*sin(#1)})
}
\draw (\TZVARxc,\TZVARyc) circle[radius=\TZVARray];
\draw \TZFONext{150}{20}--\TZFONext{110}{20}--\TZFONext{70}{20}--\TZFONext{30}{20};
\draw[color=blue, line width=0.1pt] \TZFONext{150}{20}--++\TZFONborder{75}{1};
\draw \TZFONext{150}{20}--++\TZFONborder{75}{0.69}--++\TZFONborder{40}{1.6};
\draw[color=blue, line width=0.1pt] \TZFONext{110}{20}--++\TZFONborder{60}{1.3};
\draw \TZFONext{110}{20}--++\TZFONborder{60}{1.15}--++\TZFONborder{0}{1.51};
\draw[color=blue, line width=0.1pt] \TZFONext{70}{20}--++\TZFONborder{45}{1.6};
\draw \TZFONext{70}{20}--++\TZFONborder{45}{1.41}--++\TZFONborder{-40}{1.48};
\draw [color=blue, line width=0.1pt] \TZFONext{30}{20}--++\TZFONborder{30}{1.7};
\draw \TZFONext{30}{20}--++\TZFONborder{30}{1.49};
\draw (-2.2,1.4)node[right]{\tiny successive intersection lines};
\end{scope}

\end{tikzpicture}
\end{figure}
}%
Taking the limit $\delta\tcour\to0$ we find for the intersection line the expression $\gamma(\tcour)+\RR\cdot \normal(\tcour)\wedge \normal'(\tcour)$ (in good circumstances, as in general this limit does not exists; our arguments to find the equations  \recalf{DescriptionParametriqueSurfaceReglee} and \recalf{EquationMouvementEspaceTotalBascualnt} thus are heuristic!). 
These lines are the lines of the ruled surface, envelope of the affine tangent planes $H_\tcour$, which gives us the parametric description of this developable ruled surface as
\begin{moneq}[DescriptionParametriqueSurfaceReglee]
(\tcour,u)\mapsto \gamma(\tcour) + u\cdot \normal(\tcour)\wedge \normal'(\tcour)
\mapob,
\end{moneq}
a description that (obviously) has singularities. 
Now the curve $\gamma$ also lies on the ruled surface and the fact that it is a \stress{developable} ruled surface means that we can flatten it without distortion onto a plane. 
And Persico then explains that when rolling the surface $M$ on a plane, one flattens at the same time this developable ruled surface onto this plane. 
And he concludes that parallel transport corresponds to the parallel transport in this plane in the sense of affine geometry. 

Levi-Civita uses exactly the same explanation, with as single difference that he starts with the description of the developable ruled surface, the envelope of the tangent planes along the curve, and he finishes with the idea of rolling the surface on a plane. 
Our description of the trace curve $\TraceAff\gamma\NoTraceAff{\gammah_\tfixe}$ \myquote{thus} is the result of flattening this developable ruled surface on a (fixed, affine) tangent plane, without passing through the intermediate of of the ruled surface. 
This avoids the problems that appear when the description of this surface has singularities, and at the same time it allows an easy generalisation to higher dimensions.

The arguments that lead us to the parametric description \recalf{DescriptionParametriqueSurfaceReglee} of the developable ruled surface also allows us to give a more explicit description of the isometries $\esptanbast(\tfixe,\tcour):\RR^3\to\RR^3$. 
The conditions (\refmetnaam{\labelRTS}{RTS1}) and (\refmetnaam{\labelRTS}{RTS2}) describe the movement of the rolling tangent space. 
When we look at the movement when we \myquote{roll} a tangent space $H_\tcour$ to a close by tangent space $H_{\tcour+\delta \tcour}$, it is \myquote{obvious} that this is done by a rotation $R(\tcour,\delta \tcour)$ about an axis whose direction is given by $\normal(\tcour)\wedge \normal(\tcour+\delta \tcour)$ and by an angle $\alpha$ détermined by $\cos\alpha = \inprod{\normal(\tcour)}{\normal(\tcour+\delta \tcour)}$. 
But such a rotation not only maps $H_\tcour$ (isometrically) onto $H_{\tcour+\delta \tcour}$, it also maps the whole of $\RR^3$ isometrically to itself. 
A point that at time $\tcour$ is located at $x(\tcour)\in \RR^3$ thus will be, at time $\tcour+\delta\tcour$, sent to the point $x(\tcour+\delta \tcour) = R(\tcour,\delta \tcour)\bigl(x(\tcour)\bigr)$. 
Assuming we are in the same favorable circumstances as needed in our deduction of equation \recalf{DescriptionParametriqueSurfaceReglee} for the ruled surface, we take the limit $\delta\tcour\to0$, which then tells us that the curves $x:I\to \RR^3$ must satisfy the differential equation\footnote{It is this (inhomogeneous first order linear differential) equation that is given in \cite{VanEst..1979}. The associated homogeneous equation  $v'=(\normal\wedge\normal')\wedge v$, which describes parallel transport (simply because it says that $v'$ is orthogonal to the tangent space), is exactly the equation $\extder \mathbf{R} = \pmb{\omega}\wedge\mathbf{R}$ that one finds in Levi-Civita \cite[p105]{LeviCivita..1927}.}
\begin{moneq}[EquationMouvementEspaceTotalBascualnt]
x'(\tcour) = \bigl(\normal(\tcour)\wedge \normal'(\tcour)\bigr)\wedge \bigl(x(\tcour) - \gamma(\tcour)\bigr)
\mapob.
\end{moneq}
Once we have this equation, it is easy to verify that, for $\tfixe\in I$ and $x_o\in \RR^3$ fixed, the curves $\tcour\mapsto \esptanbast(\tcour,\tfixe)(x_o)$ are the solutions of \recalf{EquationMouvementEspaceTotalBascualnt}, just as the curves $\tcour\mapsto \esptanbas(\tcour,\tfixe)(x_o)$, for $\tfixe\in I$ and $x_o\in H_\tfixe$ fixed, are the solutions of the differential equation (\refmetnaam{\labelRTS}{RTS2}).

\section{Geodesics}

A very intuitive definition of a geodesic says that it is a curve $\gamma:I\to M$ that realises \myquote{the shortest distance} between two of its points $x=\gamma(a)$ and $y=\gamma(b)$ in $M$, where the distance is the length of $\gamma$ between those two points:
\begin{moneq}
L_{\intervff ab}(\gamma) = \int_a^b \norme{\gamma'(t)}\ \extder t
\mapob.
\end{moneq}
In (every affine subspace of) $\RR^\nu$, it is a direct corollary of the triangle inequality that geodesics are (pieces of) straight lines. 
One can recognise them by the fact that their speed $\gamma'(\tcour)$ always points in the same direction, \ie, the normalised speed vector\footnote{The normalised speed vector $\gamma'(\tcour)/\norme{\gamma'(\tcour)}$ is directly linked to the curvature vector $k_\gamma(\tcour)$ given by $k_\gamma(\tcour) = \bigl(\gamma'(\tcour)/\norme{\gamma'(\tcour)}\,\bigr)'/\norme{\gamma'(\tcour)}$. Its norm $\norme{k_\gamma(\tcour)}$ is the (absolute) curvature of $\gamma$ at the point $\gamma(\tcour)$; it is the inverse of the radius of the osculating circle of $\gamma$ at this point. In $\RR^\nu$ a curve $\gamma$ thus is a geodesic if and only if its curvature vector is (identically) zero.}
\begin{moneq}
\frac{\gamma'(\tcour)}{\norme{\gamma'(\tcour)}}
\end{moneq}
is constant (independent of $\tcour$). 

For an arbitrary submanifold $M\subset \RR^\nu$ we have put \myquote{shortest distance} in quotes, as in general there is neither existence nor uniqueness. 
For non-existence it suffices to think of the plane $\RR^2$ from which a closed disc has been deleted and for non-uniqueness it suffices to think of the earth: all meridians give the same distance from north pole to south pole. 
\WelNietTikZpicture
{%
\begin{figure}[!ht]
\begin{tikzpicture}[scale=0.4, line cap=round,line join=round,x=\TZschaal, y=\TZschaal]
\draw[fill, color=lightgray] (0,0) circle (1);
\def\TZangle{20}
\draw (-3.5,0)--(-1, {2.5*sin(\TZangle)}) to[out=\TZangle, in=180] (0,1.1) to[out=0, in=180-\TZangle] (1,{2.5*sin(\TZangle)})--(3.5,0);

\def\TZangle{25}

\def\TZangle{27.5}
\draw (-3.5,0)--(-1, {2.5*sin(\TZangle)}) to[out=\TZangle, in=180] (0,1.42) to[out=0, in=180-\TZangle] (1,{2.5*sin(\TZangle)})--(3.5,0);

\def\TZangle{30}

\def\TZangle{35}
\draw (-3.5,0)--(-1, {2.5*sin(\TZangle)}) to[out=\TZangle, in=180] (0,1.75) to[out=0, in=180-\TZangle] (1,{2.5*sin(\TZangle)})--(3.5,0);

\def\TZangle{-20}
\draw (-3.5,0)--(-1, {2.5*sin(\TZangle)}) to[out=\TZangle, in=180] (0,-1.1) to[out=0, in=180-\TZangle] (1,{2.5*sin(\TZangle)})--(3.5,0);

\def\TZangle{-25}

\def\TZangle{-27.5}
\draw (-3.5,0)--(-1, {2.5*sin(\TZangle)}) to[out=\TZangle, in=180] (0,-1.42) to[out=0, in=180-\TZangle] (1,{2.5*sin(\TZangle)})--(3.5,0);

\def\TZangle{-30}

\def\TZangle{-35}
\draw (-3.5,0)--(-1, {2.5*sin(\TZangle)}) to[out=\TZangle, in=180] (0,-1.75) to[out=0, in=180-\TZangle] (1,{2.5*sin(\TZangle)})--(3.5,0);

\draw (-3.5,0)--(-1.1,0);
\draw (1.1,0)--(3.5,0);
\draw[dashed] (-1,0)--(1,0);

\draw (-3.5,0) node {$\scriptstyle\bullet$};
\draw (3.5,0) node {$\scriptstyle\bullet$};

\begin{scope}[yshift=0, xshift=10.75\TZschaal,scale=1.8]
\def\TZemme{0.7}
\def\TZhoek{30}

\draw[line width=0.3pt, color=gray] (0,0) circle (1);
\begin{scope}[param3d]

\def\TZVARx{cos(\TZphi )*sin(\TZtheta r)}
\def\TZVARy{sin(\TZphi )*sin(\TZtheta r)}
\def\TZVARz{cos(\TZtheta r)}

\def\TZphi{0}
\def\TZangleb{2.5}%
\draw[line width=0.8pt, smooth, samples=20, domain=0:\TZangleb] plot[variable=\TZtheta] ({\TZVARx}, {\TZVARy}, {\TZVARz});
\draw[dotted, line width=0.8pt, smooth, samples=20, domain=\TZangleb:pi] plot[variable=\TZtheta] ({\TZVARx}, {\TZVARy}, {\TZVARz});

\def\TZphi{80}
\def\TZangleb{2.3}
\draw[line width=0.8pt, smooth, samples=20, domain=0:\TZangleb] plot[variable=\TZtheta] ({\TZVARx}, {\TZVARy}, {\TZVARz});
\draw[dotted, line width=0.8pt, smooth, samples=20, domain=\TZangleb:pi] plot[variable=\TZtheta] ({\TZVARx}, {\TZVARy}, {\TZVARz});

\def\TZphi{-40}
\def\TZangleb{2.0}
\draw[line width=0.8pt, smooth, samples=20, domain=0:\TZangleb] plot[variable=\TZtheta] ({\TZVARx}, {\TZVARy}, {\TZVARz});
\draw[dotted, line width=0.8pt, smooth, samples=20, domain=\TZangleb:pi] plot[variable=\TZtheta] ({\TZVARx}, {\TZVARy}, {\TZVARz});

\def\TZphi{160}
\def\TZangleb{0.8}
\draw[line width=0.8pt, smooth, samples=20, domain=0:\TZangleb] plot[variable=\TZtheta] ({\TZVARx}, {\TZVARy}, {\TZVARz});
\draw[dotted, line width=0.8pt, smooth, samples=20, domain=\TZangleb:pi] plot[variable=\TZtheta] ({\TZVARx}, {\TZVARy}, {\TZVARz});

\def\TZphi{220}
\def\TZangleb{0.5}
\draw[line width=0.8pt, smooth, samples=20, domain=0:\TZangleb] plot[variable=\TZtheta] ({\TZVARx}, {\TZVARy}, {\TZVARz});
\draw[dotted, line width=0.8pt, smooth, samples=20, domain=\TZangleb:pi] plot[variable=\TZtheta] ({\TZVARx}, {\TZVARy}, {\TZVARz});

\draw (0,0,1) node {$\scriptstyle\bullet$};
\draw[color=lightgray] (0,0,-1) node {$\scriptstyle\bullet$};
\end{scope}

\end{scope}

\begin{scope}[yshift=0, xshift=20\TZschaal,scale=1.8]
\def\TZemme{0.4}
\def\TZhoek{0}

\draw[line width=0.3pt, color=gray] (0,0) circle (1);
\begin{scope}[param3d]

\def\TZVARx{cos(\TZphi )*sin(\TZtheta r)}
\def\TZVARy{sin(\TZphi )*sin(\TZtheta r)}
\def\TZVARz{cos(\TZtheta r)}

\def\TZtheta{pi/2}
\def\TZangleb{-90}\def\TZanglee{90}
\draw[line width=0.3pt, smooth, samples=20, domain=\TZangleb:\TZanglee] plot[variable=\TZphi] ({\TZVARx}, {\TZVARy}, {\TZVARz});
\draw[dotted, line width=0.3pt, smooth, samples=20, domain=\TZanglee:\TZangleb+360] plot[variable=\TZphi] ({\TZVARx}, {\TZVARy}, {\TZVARz});

\def\TZpoint{-20}
\draw ({cos(\TZpoint)},{sin(\TZpoint)},0) node {$\scriptscriptstyle\bullet$};
\def\TZpoint{20}
\draw ({cos(\TZpoint)},{sin(\TZpoint)},0) node {$\scriptscriptstyle\bullet$};

\draw[line width=1.3pt, smooth, samples=20, domain=-20:20] plot[variable=\TZphi] ({\TZVARx}, {\TZVARy}, {\TZVARz});

\end{scope}

\begin{scope}[xshift=3\TZschaal]
\draw[line width=0.3pt, color=gray] (0,0) circle (1);
\begin{scope}[param3d]

\def\TZVARx{cos(\TZphi )*sin(\TZtheta r)}
\def\TZVARy{sin(\TZphi )*sin(\TZtheta r)}
\def\TZVARz{cos(\TZtheta r)}

\def\TZtheta{pi/2}
\def\TZangleb{-90}\def\TZanglee{90}
\draw[line width=1.3pt, smooth, samples=20, domain=\TZangleb:-20] plot[variable=\TZphi] ({\TZVARx}, {\TZVARy}, {\TZVARz});
\draw[line width=1.3pt, smooth, samples=20, domain=20:\TZanglee] plot[variable=\TZphi] ({\TZVARx}, {\TZVARy}, {\TZVARz});
\draw[dotted, line width=1.3pt, smooth, samples=20, domain=\TZanglee:\TZangleb+360] plot[variable=\TZphi] ({\TZVARx}, {\TZVARy}, {\TZVARz});

\def\TZpoint{-20}
\draw ({cos(\TZpoint)},{sin(\TZpoint)},0) node {$\scriptscriptstyle\bullet$};
\def\TZpoint{20}
\draw ({cos(\TZpoint)},{sin(\TZpoint)},0) node {$\scriptscriptstyle\bullet$};

\draw[line width=0.3pt, smooth, samples=20, domain=-20:20] plot[variable=\TZphi] ({\TZVARx}, {\TZVARy}, {\TZVARz});

\end{scope}
\end{scope}

\end{scope}

\draw (0,-3) node{\footnotesize avoiding a hole};
\draw (10.75,-3) node{\footnotesize non-uniqueness};
\draw (22.7,-3) node{\footnotesize a short and a long curve};

\end{tikzpicture}
\end{figure}
}%
And even \myquote{shortest} is problematic, as a chain of \myquote{shortest} paths can lead to a \myquote{longest} path. 
For instance on the equator on the earth we can go by two ways from one point to another: up front and via the back side. 
Both are considered geodesics, but one is longer than the other. 
The official definition one adopts is that a geodesic is a curve that is stationary for the distance between two of its points: if $\gamma_\varepsilon:I\to M$ is a family of curves that all go from $x=\gamma_\varepsilon(a)$ to $y=\gamma_\varepsilon(b)$, with $\gamma_0=\gamma$, then one requires
\begin{moneq}
0 = 
\frac{\extder}{\extder \varepsilon} \bigrestricted_{\varepsilon=0} L_{\intervff ab}(\gamma_\varepsilon) 
= 
\int_a^b \fracp{\,\norme{\gamma_\varepsilon'(\tcour)}}{\varepsilon}\ \bigrestricted_{\varepsilon=0}\ \extder \tcour
\mapob.
\end{moneq}
An integration by parts gives
\begin{moneq}
\int_a^b \fracp{\,\norme{\gamma_\varepsilon'(\tcour)}}{\varepsilon}\ \bigrestricted_{\varepsilon=0}\ \extder \tcour
= - \int_a^b \inprod[4]{\frac{\partial}{\partial \varepsilon} \bigrestricted_{\varepsilon=0}\gamma_\varepsilon(\tcour)}{\frac{\extder}{\extder \tcour}\hbiggl( \frac{\gamma'(\tcour)}{\norme{\gamma'(\tcour)}}\hbiggr)} \ \extder \tcour
\mapob,
\end{moneq}
from which one deduces that this is zero for all families $\gamma_\varepsilon$ if and only if $\bigl(\gamma'(\tcour)/\norme{\gamma'(\tcour)}\bigr)'$ is orthogonal to the tangent space $T_{\gamma(\tcour)}M$, \ie, when we have
\begin{moneq}[ConditionGeodesiqueSSIDeriveeCovarianteNulle]
0 = \pi\biggl(\frac{\extder}{\extder \tcour}\hbiggl(\frac{\gamma'(\tcour)}{\norme{\gamma'(\tcour)}}\hbiggr)(\tcour)\biggr) 
\oversettext{\ \recalf{DefinitionUsuelleDeriveeCovariante}\ }\to= 
\frac{D(\gamma'/\norme{\gamma'}\,)}{\extder \tcour}(\tcour)
\mapob.
\end{moneq}
Said differently, $\gamma$ is a geodesic if and only if the covariant derivative of $\gamma'/\norme{\gamma'}$ is zero, \ie, if $\gamma'/\norme{\gamma'}$ is a field of parallel vectors along $\gamma$. 
But $\gamma'$ is a field of tangent vectors along $\gamma$ and according to \recalf{TransportparalleleDeGammahpEtGammahsec} its trace field $\TraceAff(\gamma')$ is the field $(\TraceAff\gamma)'$. 
As $\transportparallel(\tfixe,\tcour)$ is an isometry, it follows that $\norme{\gamma'} = \norme{\TraceAff\gamma'}$, which gives us the equality
\begin{moneq}
\TraceAff\hbiggl(\frac{\gamma'}{\norme{\gamma'}}\hbiggr) = \frac{\TraceAff\gamma'}{\norme{\TraceAff\gamma'}}
\mapob.
\end{moneq}
We then use \recalf{TraceDeriveeCovarianteEstDeriveeOrdinaireTrace} to obtain the equality
\begin{moneq}
\TraceAff\biggl(\frac{D(\gamma'/\norme{\gamma'}\,)}{\extder \tcour}\biggr)
=
\frac{\extder \bigl(\TraceAff(\gamma'/\norme{\gamma'}\,)\bigr)}{\extder \tcour} 
\equiv
\biggl(\frac{\TraceAff\gamma'}{\norme{\TraceAff\gamma'}}\biggr)'
\mapob.
\end{moneq}
Combined with \recalf{ConditionGeodesiqueSSIDeriveeCovarianteNulle} and the fact that $H_\tfixe$ is an affine subspace of $\RR^\nu$ in which straight lines are characterised as being curves whose normalised speed is constant, we conclude that $\gamma$ is a geodesic if and only if its trace curve $\TraceAff\gamma$ (in the affine tangent space $H_\tfixe$) is a straight line. 

\medskip

We have defined geodesics as curves $\gamma$ for which the covariant derivative of $\gamma'/\norme{\gamma'}$ is zero. 
But this is not the definition one usually finds in the literature, even though it is this condition that is invariant under an arbitrary change of parameters. 
Most authors only consider curves that are parametrised with constant speed (\ie, $\norme{\gamma'(\tcour)}$ constant), in which case the condition for being a geodesic becomes the condition that the covariant derivative of $\gamma'$ should be zero. 
Conversely, one can show that if the covariant derivative of $\gamma'$ is zero, then necessarily $\gamma$ is parametrised with constant speed. 
Defining a geodesic as a curve for which the covariant derivative of $\gamma'$ is zero simplifies the fomul\ae{} and has the additional benefit that for any $m\in M$ and any $v\in T_mM$ there exists a \stress{unique} geodesic $\gamma$ with $\gamma(0)=m$ and $\gamma'(0)=v$.

\section{By way of conclusion}

In Riemannian geometry one introduces the notion of a geodesic as a generalisation of a shortest curve (at least locally) and one introduces the notions of parallel transport of a (tangent) vector along a curve and of a field of parallel (tangent) vectors (in most cases without a real geometric interpretation). 
Imagine now a short-sighted inhabitant of a submanifold $M\subset \RR^\nu$ who thinks he is living in an affine euclidean space $H_t$ and who does not realise that, when walking, this affine space changes. 
When this inhabitant notes down his walks in his affine space (the rolling tangent space), he \stress{sees} that the shortest path between two points is a straight line and that he can translate a vector (in his affine space) from one point (of attachment) to another, which allows him to define a field of parallel vectors as one that is constant. 
The origin of this \myquote{simplified vision} lies in the fact that our inhabitant does not see the real curve $\gamma$, but only the trace curve $\TraceAff\gamma$ in his affine space $H_\tfixe$, as well as the trace vector field $\TraceAff v$ along the trace curve. 
The person inside the transparent ball who walks on the earth only sees the trace he leaves on the ground. 
Of course, having such a simplified vision of the surrounding world does not imply that our inhabitant will be unable to note that his world is not an affine space, as he will be able to observe anomalies. 
For instance that the sum of the angles in a triangle does not add up to $\pi$, or that the ratio between diameter and circumference of a circle is different from $\pi$ (see \cite{JPPGeometricon} for an amusing account of these phenomena). 
Nevertheless, for our short-sighted inhabitant geodesics really are straight lines and parallel transport is just attaching the same vector to different points.

\section{Some mathematical details}

\begin{proclaim}[ETBetTPEstIsometrie]{Lemma}
The curve $x_{s,x_o}$ (with $s\in I$ and $x_o\in H_s$) is unique if it exists. 
Moreover, the movement of the rolling tangent space  $\esptanbas(t,s):H_s\to H_t$ is an isometry:
\begin{moneq}
\forall x,y\in H_s
\quad:\quad
\norme{x-y} = \norme{\esptanbas(t,s)(x)-\esptanbas(t,s)(y)}
\mapob.
\end{moneq}

\end{proclaim}

\begin{preuve}
Let $x_{s,x_o}, y_{s,y_o}:I\to \RR^\nu$ be two such curves with $s\in I$ and $x_o,y_o\in H_s$. 
We then consider the function $f:I\to \RR$ defined by
\begin{moneq}
f(t) =\norme{x_{s,x_o}(t)-y_{s,y_o}(t)}^2
\end{moneq}
and we compute its derivative $f'(t)$:
\begin{moneq}
f'(t)
=
2\cdot \inprod{x_{s,x_o}'(t)-y_{s,y_o}'(t)}{x_{s,x_o}(t)-y_{s,y_o}(t)}
\mapob.
\end{moneq}
We then note that by \refmetnaam{\labelRTS}{RTS1} the vector $x_{s,x_o}(t)-y_{s,y_o}(t)$ belongs to $T_{\gamma(t)}M$ and that by \refmetnaam{\labelRTS}{RTS2} the vectors $x_{s,x_o}'(t)$ and $y_{s,y_o}'(t)$ are orthogonal to $T_{\gamma(t)}M$. 
It follows that $f'(t)\equiv 0$ and thus $f$ is a constant function. 
Taking $y_o=x_o$ and the initial conditions $x_{s,x_o}(s)=x_o=y_{s,x_o}(s)$ we find 
\begin{moneq}
\norme{x_{s,x_o}(t)-y_{s,x_o}(t)}^2
=
\norme{x_{s,x_o}(s)-y_{s,x_o}(s)}^2
=
\norme{x_o-x_o}^2
=0
\mapob,
\end{moneq}
hence $x_{s,x_o}=y_{s,x_o}$, proving uniqueness. 
It then suffices to change notation in the general case to find
\begin{moneq}
\norme{x_o-y_o}^2
= f(s) = f(t) =
\norme{x_{s,x_o}(t)-y_{s,y_o}(t)}^2
=
\norme{\esptanbas(t,s)(x_o)-\esptanbas(t,s)(y_o)}^2
\mapob.
\end{moneq}

\end{preuve}

\begin{proclaim}[AppIsometrieDoncAffine]{Corollary}
For fixed $s,t\in I$, the map $\esptanbas(\tcour,\tfixe):H_\tfixe\to H_\tcour$ is affine. 

\end{proclaim}

\begin{preuve}
We start with a generic computation with arbitrary $x,y,z\in H_\tfixe$: using the polarisation formula $\inprod xy = \tfrac12\cdot(\,\norme x^2+\norme y^2-\norme{x-y}^2\,)$ and writing just $E$ for $\esptanbas(\tcour,\tfixe)$, we find
\begin{align}
\shifttag{5em}
\inprod{E(x)-E(z)}{E(y)-E(z)}
\notag
\\&
=
\kern0.5em
\tfrac12\cdot\bigl(\, \norme[2]{E(x)-E(z)}^2+\norme[2]{E(y)-E(z)}^2-\norme[2]{E(x)-E(y)}^2\,\bigr)
\notag
\\&
\oversetalign{\recalt{ETBetTPEstIsometrie}}\to=
\kern0.5em
\tfrac12\cdot\bigl(\, \norme[2]{x-z}^2+\norme[2]{y-z}^2-\norme[2]{x-y}^2\,\bigr)
=
\inprod{x-z}{y-z}
\mapob.
\label{EspTanBasPreserveProduitScalaire}
\end{align}
Next we choose an orthonormal basis $f_1, \dots, f_n$ of $T_{\gamma(\tfixe)}M$ and $x_o\in H_\tfixe$ arbitrary and we define the vectors $v_i=E(x_o+f_i) - E(x_o)\in T_{\gamma(\tcour)}M$. 
According to \recalf{EspTanBasPreserveProduitScalaire} we thus have the equalities
\begin{moneq}
\inprod{v_i}{v_j} = \inprod{f_i}{f_j}
\mapob,
\end{moneq}
which shows that $v_1, \dots, v_n$ form an orthonormal basis of $T_{\gamma(\tcour)}M$. 

Because the $f_i$ form an orthonormal basis, we have for any $x\in H_\tfixe$ the equality $x-x_o=\sum_{i=1}^n \inprod{x-x_o}{f_i}\cdot f_i$. 
We thus can compute:
\begin{align*}
\shifttag{3em}
E(x)-E(x_o)
=
\sum_{i=1}^n \inprod{E(x)-E(x_o)}{v_i}\cdot v_i
\\&
=
\sum_{i=1}^n \inprod{E(x)-E(x_o)}{E(x_o+f_i)-E(x_o)}\cdot v_i
\oversettext{\ \recalf{EspTanBasPreserveProduitScalaire}\ }\to=
\sum_{i=1}^n \inprod{x-x_o}{f_i}\cdot v_i
\mapob,
\end{align*}
which can be rewritten as
\begin{moneq}
E(x) = E(x_o) + A(x-x_o)
\mapob,
\end{moneq}
where $A:T_{\gamma(\tfixe)}M \to T_{\gamma(\tcour)}M$ is the linear map defined as
\begin{moneq}
A\Bigl(\ \sum_{i=1}^n c_i\cdot f_i\ \Bigr) = \sum_{i=1}^n c_i\cdot v_i
\mapob.
\end{moneq}
This shows that $E\equiv \esptanbas(\tcour,\tfixe)$ is an affine map. 
\end{preuve}

\begin{proclaim}[ExistenceLocaleReperesMobiles]{Lemma (local existence of moving frames)}
For all $\tfixe\in I$ there exists an interval $J\subset \RR$ with $\tfixe\in J\subset I$ and $n$ maps $e_i:J\to \RR^n$, $i=1, \dots, n$ of class $C^1$ such that for all $\tcour\in J$ we have
\begin{moneq}
\{e_1(\tcour), \dots, e_n(\tcour)\}\text{ a basis of } T_{\gamma(\tcour)}M
\mapob.
\end{moneq}

\end{proclaim}

\begin{remarque}{Remark}
The family of maps $e_i:J\to \RR^\nu$ is called a \stresd{moving frame} along the curve $\gamma$; it is highly non-unique. 
A proof of its existence is a direct consequence of the existence of local charts for a submanifold.
\end{remarque}

\begin{definition}[IngredientsUtilesPourExistenceEnCoordonnees]{Some definitions useful for existence}
We start with the hypothesis that we dispose of an interval $J\subset I$ on which a moving frame $(e_1, \dots, e_n)$ is defined.
Using this moving frame, we can translate the condition (\refmetnaam{\labelRTS}{RTS1}) as the condition that there exist $n$ functions (coordinates) $\lambda_i:J\to \RR$ such that
\begin{moneq}
x(\tcour) = \gamma(\tcour) + \sum_{i=1}^n e_i(\tcour)\cdot \lambda_i(\tcour)
\mapob.
\end{moneq}
Using this expression, the condition (\refmetnaam{\labelRTS}{RTS2}) becomes the system of $n$ differential equations
\begin{moneq}
\forall 1\le j\le n
\quad:\quad
\inprod[3]{\gamma'(\tcour) + \sum_{i=1}^n e_i(\tcour)\cdot \lambda_i'(\tcour) + \sum_{i=1}^n e_i'(\tcour)\cdot \lambda_i(\tcour)}{e_j(\tcour)} = 0
\mapob.
\end{moneq}
Writing $\mfdmetric_{ji}(\tcour) = \inprod{e_j(\tcour)}{e_i(\tcour)}$, which is an invertible symmetric matrix, and $\mfdmetrichat_{ji}(\tcour)$ its inverse matrix, we can rewrite these differential equations as
\begin{moneq}
\lambda_i'(\tcour) = 
-\sum_{j=1}^n \mfdmetrichat_{ij}(\tcour)\cdot \inprod{e_j(\tcour)}{\gamma'(\tcour)} - \sum_{j,k=1}^n \mfdmetrichat_{ij}(\tcour)\cdot \inprod{e_j(\tcour)}{e_k'(\tcour)}\cdot \lambda_k(\tcour)
\mapob.
\end{moneq}
Now this is  a system of $n$ linear inhomogeneous ordinary differential equations depending on time, \ie, of the form:
\begin{moneq}
\lambda'(\tcour) = B(\tcour)\cdot \lambda(\tcour) + b(\tcour)
\mapob,
\end{moneq}
where the matrix $B(\tcour)\in \GL(n,\RR)$ and the vector $b(\tcour)\in \RR^n$ are given by
\begin{moneq}[DefvecteurPetitB]
b_i(\tcour) = -\sum_{j=1}^n \mfdmetrichat_{ij}(\tcour)\cdot \inprod{e_j(\tcour)}{\gamma'(\tcour)}
\end{moneq}
and
\begin{moneq}[DefMatriceGrandB]
B_{ik}(\tcour) = -\sum_{j=1}^n \mfdmetrichat_{ij}(\tcour)\cdot \Gamma_{jk}(\tcour)
\quad\text{with}\quad
\Gamma_{ij}(\tcour) = \inprod{e_i(\tcour)}{e_j'(\tcour)}
\mapob.
\end{moneq}
And for such equations we have the following existence and uniqueness result. 

\end{definition}

\begin{proclaim}[EquaDiffLinSOLglobalesPourGeodesiquesBasculants]{Theorem}
Let $J\subset \RR$ be an open interval and let $B:J\to \Mat(n\times n, \RR)$ and $b:J\to \RR^n$ be maps of class $C^k$, $k\ge0$ for which we consider the system of linear inhomogeneous time-dependent ordinary differential equations
\begin{moneq}[EquationDifferentielleLineaireInhomogeneGenerique]
\forall \tcour\in J
\quad:\quad
\lambda'(\tcour) = B(\tcour)\cdot \lambda(\tcour) + b(\tcour)
\mapob.
\end{moneq}
Then there exists a map $X:J\to \GL(n,\RR)$ of class $C^{k+1}$ with the following properties.\footnote{The proof of the existence of $X$ satisfying \recalf{EquaDiffLinSOLglobalesPourGeodesiquesBasculantsP1} can be found for instance in \cite[Thm. I.4.1]{L-F..1977T4}, \cite[\S XII.1]{Dieudonne..1980} or \cite[Thm. 1.6.6]{BG87}. Property \recalf{EquaDiffLinSOLglobalesPourGeodesiquesBasculantsP2} is a direct consequence.}  
\begin{enumerate}
\item \label{EquaDiffLinSOLglobalesPourGeodesiquesBasculantsP1}

The map $X$ satisfies the differential equation $X'(\tcour) = B(\tcour)\cdot X(\tcour)$.

\item\label{EquaDiffLinSOLglobalesPourGeodesiquesBasculantsP2}

When we define the maps $A:J\times J\to \GL(n,\RR)$ and $a:J\times J\to \RR^n$ by 
\begin{moneq}
A(\tcour,\tfixe) = X(\tcour)\cdot X(\tfixe)\mo
\qquad\text{and}\qquad
a(\tcour,\tfixe) = \int_{\tfixe}^{\tcour} A(\tcour,\tau)\cdot b(\tau)\ \extder \tau
\mapob,
\end{moneq}
then, for all $\tfixe\in J$ and all $\lambda_o\in \RR^n$, the map $\Lambda_{\tfixe,\lambda_o}:I\to \RR^n$
\begin{moneq}
\Lambda_{\tfixe,\lambda_o}(\tcour)= A(\tcour,\tfixe)\cdot\lambda_o+a(\tcour,\tfixe)
\end{moneq}
is the unique solution of \recalf{EquationDifferentielleLineaireInhomogeneGenerique} satisfying the initial condition $\Lambda_{\tfixe,\lambda_o}(\tfixe) = \lambda_o$. 

\end{enumerate}

\end{proclaim}

\begin{proclaim}[ExpressionExpliciteBasulementTangentEtTransportParallel]{Corollary}
For all $\tfixe\in I$ and all $x_o\in H_\tfixe$ there exists a unique curve $x:I\to \RR^\nu$ satisfying \refmetnaam{\labelRTS}{RTS1}, \refmetnaam{\labelRTS}{RTS2} and the initial condition $x(\tfixe) = x_o$.

\end{proclaim}

\begin{preuve}
For all $\tfixe \in I$ there exists an interval $J$ such that $\tfixe\in J\subset I$ and on which there exists a moving frame \recalt{ExistenceLocaleReperesMobiles}. 
According to \recalt{IngredientsUtilesPourExistenceEnCoordonnees} and \recalt{EquaDiffLinSOLglobalesPourGeodesiquesBasculants} there exists a unique solution $x_{\tfixe,x_o}:J\to \RR^\nu$ satisfying (\refmetnaam{\labelRTS}{RTS1}) and (\refmetnaam{\labelRTS}{RTS2}) as well as $x_{\tfixe,x_o}(\tfixe) = x_o$. 
Now suppose $J=\intervoo ab\subset I$ is the largest (open) interval on which such a solution exists. 
If $J\neq I$, then either $a\in I$ or $b\in I$, so suppose it is $b\in I$. 
Again by \recalt{ExistenceLocaleReperesMobiles} there exists an open interval  $J'$ with $b\in J'\subset I$ on which a moving frame exists. 
We choose arbitrarily $\tfixe'\in J \cap J'$ and we define $y_o=x_{\tfixe,x_o}(\tfixe')\in H_{\tfixe'}$.
By the previous argument there exists a unique solution $x_{\tfixe',y_o}:  J'\to \RR^\nu$ satisfying (\refmetnaam{\labelRTS}{RTS1}) and (\refmetnaam{\labelRTS}{RTS2}) as well as $x_{\tfixe',y_o}(\tfixe') = y_o$. 
But the curve $x_{\tfixe,x_o}$ satisfies these same requirements (on $J$). 
By uniqueness in \recalt{ETBetTPEstIsometrie} it follows that the two curves $x_{\tfixe',y_o}$ and $x_{\tfixe,x_o}$ coincide on the intersection $J\cap J'$, which allows us to extend the definition of $x_{\tfixe,x_o}$ to $J\cup J'$, contradicting maximality of $J$. 
This shows that we can not have $b\in I$. 
A similar argument shows that we can not have $a\in I$ and thus $J=I$. 
\end{preuve}

\begin{remarque}[RemarqueExpressionLineaireMouvEspTanBasc]{Remark}
If there exists a moving frame on an interval $J\subset I$, then the movement of the rolling tangent space $\esptanbas(\tcour,\tfixe)$ and the parallel transport map $\transportparallel(\tcour,\tfixe)$ are given by
\begin{align*}
\esptanbas(\tcour,\tfixe)\biggl(\gamma(\tfixe) + \sum_{i=1}^n  e_i(\tfixe)\cdot\lambda_i\biggr) 
&
=
\gamma(\tcour) + \sum_{i=1}^n e_i(\tcour)\cdot\biggl(a_i(\tcour,\tfixe)+\sum_{j=1}^n A_{ij}(\tcour,\tfixe)\cdot \lambda_j\biggr)
\\[2\jot]
\transportparallel(\tcour,\tfixe)\biggl(\, \sum_{i=1}^n  e_i(\tfixe)\cdot v_i\biggr)
&
=
\sum_{i=1}^n e_i(\tcour)\cdot\biggl(\,\sum_{j=1}^n A_{ij}(\tcour,\tfixe)\cdot v_j\biggr)
\mapob.
\end{align*}
These expressions confirm that $\esptanbas(\tcour,\tfixe)$ is an affine map with $\transportparallel(\tcour,\tfixe)$ its linear part. 
\end{remarque}

\begin{proclaim}[DeriverEspTanBasstfsDonneTransportParallelPartieParallele]{Lemma}
Let $\tfixe\in I$ be fixed, let $f:I\to H_\tfixe$ be a map of class $C^1$ and let $g:I\to \RR^\nu$ be defined as $g(\tcour) = \esptanbas(\tcour,\tfixe)\bigl(f(\tcour)\bigr)$. 
Then we have $\transportparallel(\tcour,\tfixe)\bigl(f'(\tcour)\bigr) = \pi\bigl(g'(\tcour)\bigr)$ (with $\pi$ the orthogonal projection on $T_{\gamma(\tfixe)}M$).

\end{proclaim}

\begin{preuve}
We start by introducing the function $G:I\times I\to H_\tfixe$ as
\begin{moneq}
G(u,v) = \esptanbas(u,\tfixe)\bigl(f(v)\bigr)
\end{moneq}
and we note that we have $g(\tcour) = G(\tcour,\tcour)$. 
Hence we have
\begin{moneq}
g'(\tcour) = (\partial_1G)(\tcour,\tcour) + (\partial_2G)(\tcour,\tcour)
\mapob.
\end{moneq}
By definition of the rolling tangent space, the vector
\begin{moneq}
(\partial_1G)(u,v)
=
\partial_u\esptanbas(u,\tfixe)\bigl(f(v)\bigr)
\end{moneq}
is orthogonal to $T_{\gamma(u)}M$ and because $\esptanbas(\tfixe,\tcour)$ is affine, we have:
\begin{align*}
(\partial_2G)(u,v)
&
=
\lim_{h\to 0} \frac{\esptanbas(u,\tfixe)\bigl(f(v+h)\bigr) - \esptanbas(u,\tfixe)\bigl(f(v)\bigr)}{h}
\\&
=
\lim_{h\to 0} \transportparallel(u,\tfixe)\biggl(\frac{f(v+h)-f(v)}{h} \biggr)
=
\transportparallel(u,\tfixe)\bigl(f'(v)\bigr)
\mapob.
\end{align*}
Now $\transportparallel(\tcour,\tfixe)\bigl(f'(\tcour)\bigr)$ belongs to $T_{\gamma(\tcour)}M$, hence $\pi\bigl(g'(\tcour)\bigr) = (\partial_2G)(\tcour,\tcour)$. 
\end{preuve}

\begin{proclaim}[TransportparalleleDeGammahpEtGammahsec]{Proposition}
For fixed $\tfixe\in I$ the trace curve $\TraceAff\gamma\NoTraceAff{\gammah_\tfixe}:I\to H_\tfixe$ and the trace vector field $\TraceAff v\NoTraceAff{\vh_\tfixe}:I\to T_{\gamma(\tfixe)}M$ satisfy the equalities
\begin{moneq}
(\TraceAff\gamma)'(\tcour) = \bigl(\TraceAff(\gamma')\bigr)(\tcour)
\quad\text{and}\quad
\transportparallel(\tcour,\tfixe)\bigl((\TraceAff v)'(\tcour)\bigr) = \pi\bigl(v'(\tcour)\bigr)
\mapob.
\end{moneq}

\end{proclaim}

\begin{preuve}
Using \recalf{LoiDeGroupeRTB} we rewrite the definition of $(\TraceAff\gamma\NoTraceAff{\gammah_\tfixe})(\tcour) = \esptanbas(\tfixe,\tcour)\bigl(\gamma(\tcour)\bigr)$ as
\begin{moneq}
\gamma(\tcour) = \esptanbas(\tcour,\tfixe)\bigl((\TraceAff\gamma\NoTraceAff{\gammah_\tfixe})(\tcour)\bigr)
\end{moneq}
and we apply \recalt{DeriverEspTanBasstfsDonneTransportParallelPartieParallele} to obtain $\transportparallel(\tcour,\tfixe)\bigl((\TraceAff\gamma\NoTraceAff{\gammah_\tfixe})'(\tcour)\bigr) = \pi\bigl(\gamma'(\tcour)\bigr) = \gamma'(\tcour)$ (because $\gamma'(\tcour)\in T_{\gamma(\tcour)}M$), 
and thus
\begin{moneq}
\bigl(\TraceAff(\gamma')\bigr)(\tcour)
\oversettext{\ \recalf{DefinitionChampTraceLeLongCourbe}\ }\to=
\transportparallel(\tfixe,\tcour)\bigl(\gamma'(\tcour)\bigr)
=
\transportparallel(\tfixe,\tcour)\transportparallel(\tcour,\tfixe)\bigl((\TraceAff\gamma)'(\tcour)\bigr)
=
(\TraceAff\gamma)'(\tcour)
\mapob,
\end{moneq}
which is the first equality. 

For the vector field $\TraceAff v\NoTraceAff{\vh_\tfixe}$ we rewrite its definition as $v(\tcour) = \transportparallel(\tcour,\tfixe)\bigl((\TraceAff v\NoTraceAff{\vh_\tfixe})(\tcour)\bigr)$ and we apply \recalf{TransportParallelParDifferenceEspTanBas} with $x_o\in H_\tfixe$ arbitrary to obtain
\begin{moneq}
v(\tcour) = \esptanbas(\tcour,\tfixe)\bigl(x_o+(\TraceAff v\NoTraceAff{\vh_\tfixe})(\tcour)\bigr) - \esptanbas(\tcour,\tfixe)(x_o)
\mapob.
\end{moneq}
Applying \recalt{DeriverEspTanBasstfsDonneTransportParallelPartieParallele} a second time (twice, with $f_1(\tcour) = x_o$ and $f_2(\tcour)=x_o+(\TraceAff v\NoTraceAff{\vh_\tfixe})(\tcour)$), we get
\begin{align*}
\pi\bigl(v'(\tcour)\bigr) = \transportparallel(\tcour,\tfixe)\bigl((\TraceAff v\NoTraceAff{\vh_\tfixe})'(\tcour)\bigr) - \transportparallel(\tcour,\tfixe)(0) = \transportparallel(\tcour,\tfixe)\bigl((\TraceAff v\NoTraceAff{\vh_\tfixe})'(\tcour)\bigr)
\mapob,
\end{align*}
which is the second equality. 
\end{preuve}

\begin{remarque}{Remark for the connoisseurs}
In local coordinates on $M$, the equation \recalf{ConditionGeodesiqueSSIDeriveeCovarianteNulle} becomes the system of equations
\begin{align}
\shifttag{5em}
c_j''(\tcour) + \sum_{k,\ell=1}^n \Gamma^j{}_{\ell k}\bigl(c(\tcour)\bigr)\cdot c_\ell'(\tcour)\cdot  c_k'(\tcour)
\notag
\\
&=\ 
\frac{\ \displaystyle\sum_{i,j=1}^n c_i'(\tcour)\cdot \mfdmetric_{ij}\bigl(c(\tcour)\bigr)\cdot \Bigl( c_j''(\tcour) + \sum_{k,\ell=1}^n \Gamma^j{}_{\ell k}\bigl(c(\tcour)\bigr)\cdot c_\ell'(\tcour)\cdot  c_k'(\tcour) \Bigr)\ }{\displaystyle\sum_{i,j=1}^n c'_i(\tcour)\cdot \mfdmetric_{ij}\bigl(c(\tcour)\bigr)\cdot c_j'(\tcour)}
\mapob,
\label{ExpressionEquationsGeodesiquesCoordonneesLocales}
\end{align}
where the $c_i$ are the local coordinates of the curve $\gamma$, where $\mfdmetric_{ij}$ denotes the matrix of the metric (on $M$, induced by the standard metric on $\RR^\nu$) and where $\Gamma^j{}_{\ell k}$ denote the Christoffel symbols associated to this metric. 
The left hand side of \recalf{ExpressionEquationsGeodesiquesCoordonneesLocales} gives the covariant derivative of $\gamma'$ in terms of the local coordinates. 
\end{remarque}

\section*{Acknowledgements}
We are indebted to R.~Tazzioli for her invaluable help with the historical details and references. 
During the preparation of this manuscript the first author was supported by a scholarship of the IKS graduate programme of the University of Lille. 
The second author acknowledges the support of the CDP C2EMPI, as well as the French State under the France-2030 programme, the University of Lille, the Initiative of Excellence of the University of Lille, the European Metropolis of Lille for their funding and support of the R-CDP-24-004-C2EMPI project.

\providecommand{\bysame}{\leavevmode\hbox to3em{\hrulefill}\thinspace}
\providecommand{\MR}{\relax\ifhmode\unskip\space\fi MR }
\providecommand{\MRhref}[2]{%
  \href{http://www.ams.org/mathscinet-getitem?mr=#1}{#2}
}
\providecommand{\href}[2]{#2}

\end{document}

\end